\RequirePackage{fix-cm}

\documentclass[a4,12pt,twoside,openright]{article}
\usepackage{graphicx}
\usepackage[a4paper, top=3cm, bottom=3cm, left=2cm, right=2cm]{geometry}
\usepackage{bbm, amssymb, amsmath, amsfonts, fancyhdr} 
\usepackage[english]{babel}
\usepackage{array}
\usepackage{oldgerm} 
\usepackage[T1]{fontenc}
\usepackage[thmmarks]{ntheorem}
\usepackage{pst-all}
\usepackage[crop=off]{auto-pst-pdf}

\newcommand{\N}{\mathbbm N}		
\newcommand{\R}{\mathbbm R}		

\newcommand{\FT}{\textfrak F}		
\newcommand{\s}{\mathcal S}			
\newcommand{\landau}{\mathcal O}
\newcommand{\inv}{^{-1}}

\newcommand{\diam}{\mathop{diam}\nolimits}
\newcommand{\supp}{\mathop{supp}}
\newcommand{\dd}{~\mathrm{d}}
\newcommand{\minus}{\!\textrm{-}}

\newcommand{\x}{\!\times\!}
\newcommand{\halb}{\durch{2}}
\newcommand{\e}{\varepsilon}

\newcommand{\zweik}{(2^k\minus1)}\newcommand{\zweil}{(2^l\minus1)}
\newcommand{\mhalb}{\tfrac{m}{2}}
\newcommand{\krange}{0,\ldots,|\log_2\gamma|-1}
\newcommand{\zrange}{\durch{\beta},\ldots,\tfrac{2}{\beta}-1}
\newcommand{\Gammaplus}{\Gamma_{\alpha}^\delta+\Gamma_{\mu}^\delta}

\newcommand{\durch}[1]{\ensuremath{\frac{1}{#1}}}

\renewcommand{\labelenumi}{(\roman{enumi})}

\usepackage[thmmarks]{ntheorem}

\theoremstyle{plain}

\newtheorem{stepnr}{Step}

\theoremstyle{break} 
\theoremheaderfont{\scshape}
\newtheorem{satznr}{Theorem}
\newtheorem{lemnr}[satznr]{Lemma}
\newtheorem{kornr}[satznr]{Corollary}
\newtheorem{bemnr}[satznr]{Remark}

\theoremstyle{nonumberbreak} 
\newtheorem{defi}{Definition}

\newtheorem{satz}{Theorem}[section] 

\newtheorem{bem}{Remark} 

\theoremstyle{nonumberplain}
\theoremsymbol{$\Box$} 
\theorembodyfont{\upshape}
\theoremheaderfont{}
\newtheorem{bew}{Proof:\,\,}				

\newcounter{abb}

\begin{document}

\title{A sharp $L^p$-$L^q-$Fourier restriction theorem for a conical surface of finite type
}


\author{Stefan Buschenhenke 
}




\maketitle

\begin{abstract}
In this paper, we study the Fourier restriction problem for certain conical surfaces where the sections are curves of finite type.
We obtain the sharp $L^p-L^q$-range for this surfaces.
Our methods rely on a variation of the affine arclength measure.
\end{abstract}

\renewcommand\contentsname{Contents}
\tableofcontents
\thispagestyle {empty}  
\newpage



\section{Introduction and first considerations}

\subsection{Introduction}
Let $S$ be a compact hypersurface in $\R^n$ (or more generally a smooth submanifold) with surface measure $\sigma$. We say that the $L^p(\R^n)$-$L^q(S)$ \textit{Restriction estimate} holds if there exists a constant $C>0$ such that
\begin{equation}
  \left( \int\limits_S |\hat f(\xi)|^q \dd\sigma(\xi) \right)^{\tfrac{1}{q}} \leq C\cdot \|f\|_p 
\end{equation}
for every Schwartz function $f$, where
\begin{align}
	\hat f (x) = \int_{\R^n} f(\xi) e^{-ix\cdot\xi}\dd\xi,\qquad g\in L^1(\R^n),x\in\R^n,
\end{align}
denotes the Fourier transform of $f\in\s(\R^n)$.


In a very seminal achievement, Stein and Tomas proved that the restriction operator of the unit sphere $R:L^{p}(\R^n)\to L^{2}(S^{n-1})$ is bounded if and only if $1\leq p \leq\frac{2n+2}{n+3}$ [St]. The crucial property in the proof is
the non-vanishing Gaussian curvature of the sphere.\\
Concerning $(p,q)$-estimates, it is conjectured that $R:L^p(\R^n)\to L^q(S^{n-1})$ is bounded if $p'>\frac{2n}{n-1}$ and $\frac{1}{q}\geq\frac{n+1}{(n-1)p'}$. In dimension $n=2$, the conjecture was confirmed by Zygmund in 1974 [Z]. However, in higher dimensions, the conjecture is still open, despite partial results, namely by Tao [T] and most recently by Bourgain and Guth [BG].\\

The surface we will consider is a surface of so-called finite type, where the tangent plane has finite order of contact. This means - describing the surface locally as a graph of a function - that for all directions, the partial derivatives of a certain order are not vanishing. However, the second order partial derivatives of the function (and therefore the corresponding principle curvature) is allowed to vanish at some points. 
 A simple example is the curve $x^m$, $m\geq3$.
The restriction estimates for curves are known due to subsequent work of Sj\"olin [Sj], Ruiz [R] and Barcelo [B2]. But we can construct a cone-like surface from the curve, like the "classical" cone arises from a circle (which has non-vanishing curvature). A sharp restriction theorem for the "classical" cone in $\R^3$ was proven by Barcelo [Ba1].
It should be mentioned that there was earlier work on this problem by Barcelo [B2], obtaining partial results.

Our approach involves a certain weighted estimate, which can be considered as a variation of the affine arclength measure. In Fourier restriction theory of curves, the affine arclength measure has turned out to be a very strong tool, see for instance [DM].

\subsection{The main theorem and necessary conditions}
Let $\gamma$ be a compact curve of finite type, i.e. for all $p\in\gamma$ exists a local parametrisation $\Phi$, $\Phi(x)=p$. The minimal $m\in\N_{\geq2}$ such that $\Phi$ is $m-$times differentiable at $x$ and $\Phi^{(m)}(x)\neq0$ is called the type of $\gamma$ at $p$. Let $M$ be the maximal type of $\gamma$.

Our main result, given in terms of Lorentzspaces $L^{p,r}$, is the following:

\begin{satznr}\label{thm0}
Let $1\leq p,q,r\leq\infty$.
We consider the patch of a conical surface \\
$\Gamma=\{(\xi,z)\in\R^2\x\R|\,1\leq z\leq2, ~\frac{\xi}{z}\in\gamma \}$ with surface measure $\sigma$.
Then the Fourier restriction inequality
\begin{align}\label{restsatz0}
\|\hat{f} |_\Gamma\|_{L^{q,r}(\Gamma,\sigma)} \leq C \|f\|_{L^{p,r}(\R^3)} \qquad \forall f\in \s(\R^3),
\end{align}
holds if $1\leq p<\tfrac{M+1}{M}$ and $\durch{q}\geq\tfrac{M+1}{p'}$.\\
If furthermore $p\leq q$ or $\durch{q}>\tfrac{M+1}{p'}$, then we have
\begin{align}\label{restsatz01}
\|\hat{f} |_\Gamma\|_{L^{q}(\Gamma,\sigma)} \leq C \|f\|_{L^{p}(\R^3)} \qquad \forall f\in \s(\R^3).
\end{align}
\end{satznr}

The theorem is sharp in the sense that there exists surfaces where \eqref{restsatz0} is not valid if the conditions on $p'$ and $q$ are violated. That the condition $\frac{1}{q}\geq\frac{M+1}{p'}$ is necessary can be seen by a classical Knapp-type example. \\
Notice further that for curves of finite type, it is known that the strong $L^p\to L^q$-estimate fails for $p>q$ and $\durch{q}=\tfrac{M+1}{p'}$ [So]. A corresponding argument for the cone can be found in Chapter \ref{kegelnotwendig}. At the endpoint $p'=m+1$, a weak-type estimate might still hold, but this is beyond our methods.\\

Observe that \eqref{restsatz0} implies \eqref{restsatz01} for the described values of $p$ and $q$. For $p\leq q$, this is a consequence of the Marcinkiewicz interpolation theorem, whereas for $\durch{q}>\tfrac{M+1}{p'}$, we use the fact that for $q<\tilde q$ we have
\begin{align}
	\|\hat{f} |_\Gamma\|_{L^{q,1}(\Gamma,\sigma)}\leq \|\hat{f} |_\Gamma\|_{L^{\tilde q,\infty}(\Gamma,\sigma)}
\end{align}
 since $\sigma(\Gamma)<\infty$.\\


\section{Reduction of the problem}

Localising to points of vanishing curvature, and describing the surface locally as a graph, it remains to discuss the following problem:
\begin{defi}Let $m\in\{3,4,\ldots\}$ and let $\Phi:[0,1]\to\R$ satisfy the following conditions:
\begin{enumerate}\renewcommand{\labelenumi}{(\roman{enumi})}
	\item $\exists \chi \in C^2[0,1],\chi>0: \Phi(x)=x^m \chi(x)$, \\
				$\exists \chi^1 \in C^1[0,1],\chi^1>0: \Phi'(x)=x^{m-1} \chi^1(x)$, \\
				$\exists \chi^2 \in C [0,1],\chi^2>0: \Phi''(x)=x^{m-2} \chi^2(x)$,
	\item $\Phi^{(k)}(x)>0$ for $x>0$, $k\leq m$; especially, $\Phi$ and $\Phi'$ are convex.
\end{enumerate}
We define the generalised cone (more exactly, a section of a cone) 
$$\Gamma=\{(x,y,z)\in\R^3|\,0\leq x\leq 1,~1\leq z\leq 2,~\Phi\left(\frac{x}{z}\right)=\frac{y}{z}\}.$$
The associated surface measure will be denoted by $\sigma$. 
\end{defi}

\begin{satznr}\label{thm1}
For $1\leq p<\tfrac{m+1}{m}$ and $\durch{q}\geq\tfrac{m+1}{p'}$ holds 
\begin{align} \label{restsatz} 
\|\hat{f} |_\Gamma\|_{L^{q}(\Gamma,\sigma)} \leq C \|f\|_{L^{p}(\R^3)} \qquad \forall f\in \s(\R^3).
\end{align}
\end{satznr}


\subsection{Decomposition}\label{zerlegungteil}
The critical part of $\Gamma$ is the line $x=0$, where the curvature vanishes. To take this fact into account, we decompose $\Gamma$ into dyadic pieces, becoming smaller near the line $x=0$.
By rescaling to the case of (almost) constant curvature, we would be able to make use of already known estimates. Unfortunately, summation is only possible if $\durch{q}>\tfrac{m+1}{p'}$.
In the limit case $\durch{q}=\tfrac{m+1}{p'}$, we need to deal with the problem as a whole. For this, we will further decompose each dyadic piece, depending on the curvature.


Let $\Gamma^\delta=\{(x,y,z)\in\R^3|\,0\leq x< 1,~1\leq z\leq 2,~\Phi\left(\frac{x}{z}\right)\leq\frac{y}{z}\leq\Phi\left(\frac{x}{z}\right)+\delta\}$ be the thickening of $\Gamma$ by $\delta>0$ (we drop the points with $x=1$ for technical purposes).
Moving on the $x-$axis from the origin by length $\sqrt[m]{\delta}$ corresponds to $\Phi$ changing by $\delta$. In other words, this part of $\Gamma$ is contained in a box of width $\delta$.
We thus define $\gamma=\sqrt[m]{\delta}$ and
$$\Gamma_k^{\delta}=\{(x,y,z)\in\Gamma^\delta|\zweik\gamma\leq x< (2^{k+1}\minus1)\gamma\},\qquad k=0,\ldots,\left\lceil \log_2\durch\gamma\right\rceil-1.$$
Now how to determine the finer decomposition? We change coordinates, or respectively $\Phi$ by affine transformation into
$$\Phi^k(x)=\Phi(x+\zweik\gamma)-\Phi((2^k\minus1)\gamma)-x\Phi'(\zweik\gamma),$$ such that $$\Phi^k(0)=0=(\Phi^k)'(0).$$
According to Taylor, we get
$$\Phi^k(x)\approx\frac{m(m-1)}{2}(\zweik\gamma)^{m-2}x^2+\landau(x^3).$$
On which distance $\gamma_k$ from the (new) origin does $\Phi$ varies at most $\delta$?
We demand ($k\neq0$)
\begin{align*}
	\delta=&\Phi^k(\gamma_k)\approx \frac{m(m-1)}{2}(\zweik\gamma)^{m-2}\gamma_k^2, \\
	\text{i.e.}\qquad \gamma_k\approx& \sqrt{\delta (\zweik\gamma)^{2-m}} 
							\approx 2^{k(1-\tfrac{m}{2})}\gamma.
\end{align*}
Concerning the z-coordinate, we decompose equidistantly with width $\beta$, where we require $\beta\lesssim\delta$. This ensures that the projection of such a set in $x$-$y-$space does not appear to different to a intersection parallel to $x$-$y-$space. To choose $\beta=\delta$ would be appropriate and you might assume this. Nevertheless, we will distinct these two quantities to be aware how each of them effects our computations. We will see that all the $\beta$'s cancels at the end of the proof, reflecting the fact that there is no impact from the $z$-direction.
We obtain a decomposition of $\Gamma^\delta$ as follows:

\begin{defi}[Decomposition]
Let $\delta>0$, $\gamma>0$ and $\gamma^m=\delta$, let $\beta<\delta$ with $\durch{\beta}\in\N$. For $k=0,\ldots,\left\lceil \log_2\durch\gamma\right\rceil-1$, $j\in I_k=\{0,\ldots,2^{k\mhalb}-1\}$ and $n=\durch{\beta},\ldots\frac{2}{\beta}-1$ define
$\gamma_k=2^{k(1-\mhalb)}\gamma$, $x_{kj}=\zweik\gamma+j\gamma_k$, $x_k=x_{k,0}$ and
\begin{align}
	\Gamma_{kjn}^\delta = \{(x,y,z)\in\R^3|\, n\beta\leq z\leq(n+1)\beta,~
	\Phi\left(\tfrac{x}{z}\right) \leq \tfrac{y}{z}\leq \Phi\left(\tfrac{x}{z}\right)+\delta,~
  x_{kj}\leq x< x_{k,j+1}\}.
\end{align}
Furthermore let $\phi_{kjn}$ be a bump function adapted to $\Gamma_{kjn}$. To be more precise,
if $\eta\in C_0^\infty(\R)$, $\chi_{[-1,1]}\leq\eta\leq\chi_{[-\tfrac{5}{4},\tfrac{5}{4}]}$, let
$$\phi_{kjn}(x,y,z)=\eta\left(\frac{x-x_{kj}}{\gamma_k}\right)
\eta\left(\frac{y-z\Phi(x/z)}{\delta}\right)\eta\left(\frac{z-n\beta}{\beta}\right).$$
This means that $\phi_{kjn}$ is to some extend supported in an anisotropic thickening of $\Gamma_{kjn}^\delta$, precisely in the set
\begin{align*}
	\tilde\Gamma_{kjn}^\delta =& \{(x,y,z)\in\R^3|\, |z-n\beta|\leq\tfrac{5}{4}\beta,~
	|\tfrac{y}{z}-\Phi\left(\tfrac{x}{z}\right)|\leq\tfrac{5}{4}\delta,~
  |x-x_{kj}|\leq \tfrac{5}{4}\gamma_k\} \\
  \subset& \bigcup\limits_{u,v,w\in\{-1,0,1\}} (u\gamma_k,v\delta,w\beta)+\Gamma_{kjn}^\delta.	
\end{align*}
For the further proceeding, we will always denote by $\alpha$ the triple $(k,j,n)$, and, if required, by $\mu$ a second triple $(l,i,p)$.\\
For simplification, we write 
$\sum\limits_\alpha$ instead of $\sum\limits_{k=0}^{|\log_2 \gamma|-1}\sum\limits_{j=0}^{2^{k\frac{m}{2}}-1} \sum\limits_{n=\durch{\beta}}^{\frac{2}{\beta}-1}$.
\end{defi}

\stepcounter{abb}
\ifx\JPicScale\undefined\def\JPicScale{0.5}\fi
\psset{unit=\JPicScale mm}
\psset{linewidth=0.3,dotsep=1,hatchwidth=0.3,hatchsep=1.5,shadowsize=1,dimen=middle}
\psset{dotsize=0.7 2.5,dotscale=1 1,fillcolor=black}
\psset{arrowsize=1 2,arrowlength=1,arrowinset=0.25,tbarsize=0.7 5,bracketlength=0.15,rbracketlength=0.15}
\begin{pspicture}(0,0)(238.42,186.37)
\psline[arrowscale=2 1]{<-}(195.26,42.66)(55,42.66)
\psline[arrowscale=2 1]{<-}(73.09,60.51)(33.42,21.37)
\newrgbcolor{userLineColour}{0.4 0.4 0.4}
\psline[linecolor=userLineColour](238.42,117.18)(157.5,37.34)
\newrgbcolor{userLineColour}{0.4 0.4 0.4}
\psline[linecolor=userLineColour](103.97,58.77)(82.39,37.48)
\newrgbcolor{userLineColour}{0.4 0.4 0.4}
\psline[linecolor=userLineColour](146.71,58.63)(125.13,37.34)
\newrgbcolor{userLineColour}{0.4 0.4 0.4}
\psline[linecolor=userLineColour](81.97,58.63)(60.39,37.34)
\rput(38.81,32.02){y}
\rput(187.17,40){x}
\psline[arrowscale=2 1]{<-}(55,181.05)(55,42.66)
\psline(55,149.11)(55,74.6)
\psline(195.26,106.53)(200.66,111.85)
\psline(195.26,106.53)(173.68,159.76)
\psline(179.08,165.08)(200.66,111.85)
\psline[linestyle=dotted,dotsep=3](173.68,159.76)(173.68,53.31)
\psline[linestyle=dotted,dotsep=3](195.26,106.53)(195.26,74.6)
\psline[linestyle=dotted,dotsep=3](162.89,149.11)(162.89,42.66)
\psline[linestyle=dotted,dotsep=3](162.89,149.11)(55,149.11)
\newrgbcolor{userLineColour}{0.4 0.4 0.4}
\psline[linecolor=userLineColour](173.68,159.76)(162.89,149.11)
\newrgbcolor{userLineColour}{0.4 0.4 0.4}
\psline[linecolor=userLineColour](200.66,186.37)(179.08,165.08)
\newrgbcolor{userLineColour}{0.4 0.4 0.4}
\psline[linecolor=userLineColour](238.42,149.11)(162.89,74.6)
\psline[linestyle=dotted,dotsep=3](200.66,111.85)(200.66,79.92)
\psbezier(55,99.67)(179.65,99.39)(188.44,124.32)(188.44,124.32)
\psline(181.2,142.53)(186.6,147.86)
\psline(188.73,123.9)(194.13,129.23)
\psline(65.79,149.11)(66.64,74.72)
\psline(88.79,75.71)(87.51,149.39)
\psline(124.99,151.35)(132.09,81.59)
\psline[linestyle=dashed,dash=1 1](109.09,77.53)(106.25,150.23)
\psline[linestyle=dashed,dash=1 1](96.73,149.8)(98.87,76.69)
\psline[linestyle=dashed,dash=1 1](115.19,150.64)(120.16,79.07)
\psline[linestyle=dotted,dotsep=3](132.09,81.41)(132.09,43.99)
\psline[linestyle=dotted,dotsep=3](88.92,75.7)(88.92,41.55)
\rput(49.6,165.08){z}
\psbezier(60.39,154.43)(162.89,154.43)(179.08,165.08)(179.08,165.08)
\psline(173.68,159.76)(179.08,165.08)
\psline(60.39,154.43)(55,149.11)
\psbezier(55,149.11)(157.5,149.11)(173.68,159.76)(173.68,159.76)
\pscustom[]{\psbezier(55,75)(172.63,75)(195.26,106.97)(195.26,106.97)
\psbezier(195.26,106.97)(195.13,106.44)(195.26,106.97)
}
\pscustom[]{\psbezier(55,125)(152.63,124.87)(180.79,142.37)(180.79,142.37)
\psline(180.79,142.37)(180.79,142.5)
}
\rput(110,10){\textbf{Figure \thesection.\theabb:} Decomposition}
\end{pspicture}


\subsection{The heart of the problem}

Of essential impact is the following theorem. It is a weighted, discrete version of the adjoint restriction estimate and, as we will see, already implies the restriction theorem.

\begin{satznr}\label{superlemma}
Let $p'>m+1$, $\frac{3}{p'}<\tfrac{1}{q}< \tfrac{1}{2}+\tfrac{1}{p'}$. Then 
\begin{align}\label{arclength}
		\|\sum_{\alpha} a_\alpha \hat\phi_\alpha\|_{p'}	\lesssim&  
		\delta^{\durch{q}}\|\sum_{\alpha=(k,j,n)} 
							a_\alpha (\zweik\gamma)^{\durch{q}-\frac{m+1}{p'}} \phi_\alpha\|_{q'}.
\end{align}
\end{satznr}


~\\
~\\

\stepcounter{abb}
\def\JPicScale{0.46}
\psset{unit=\JPicScale mm}
\psset{linewidth=0.3,dotsep=1,hatchwidth=0.3,hatchsep=1.5,shadowsize=1,dimen=middle}
\psset{dotsize=0.7 2.5,dotscale=1 1,fillcolor=black}
\psset{arrowsize=1 2,arrowlength=1,arrowinset=0.25,tbarsize=0.7 5,bracketlength=0.15,rbracketlength=0.15}
\begin{pspicture}(0,0)(182.21,161.97)
\newrgbcolor{userLineColour}{0.8 1 1}
\newrgbcolor{userFillColour}{0.2 0.2 1}
\pspolygon[linewidth=0.35,linecolor=userLineColour,linestyle=none,fillcolor=userFillColour,fillstyle=vlines](24.19,37.07)
(155.49,116)
(23.95,89.61)
(23.95,89.85)(24.19,37.07)
\psline[linewidth=0.4,arrowsize=1.3 2]{->}(24.19,36.99)(182.21,36.99)
\psline[linewidth=0.4](24.19,142.3)(159.16,142.3)
\psline[linewidth=0.4,linestyle=dotted](155.87,142.3)(155.87,36.99)
\psline[linewidth=0.4](22.54,89.65)(25.83,89.65)
\psline[linewidth=0.35](22.54,63.32)(24.19,63.32)
\psline[linewidth=0.35](22.54,115.97)(24.19,115.97)
\psline[linewidth=0.35](71.92,142.3)(75.21,148.88)
\rput(17.96,157.88){1/q}
\rput(181.41,29.1){1/p'}
\rput(155.01,28.04){1/4}
\rput(12.87,90.43){1/2}
\rput(12.87,116.76){3/4}
\rput(12.87,64.11){1/4}
\rput(14.51,36.47){0}
\psline[linewidth=0.35,linestyle=dotted](50.52,94.91)(50.52,36.99)
\psline[linewidth=0.35,linestyle=dotted](71.92,142.3)(71.92,36.99)
\rput(80.83,27.81){1/(m+1)}
\rput(46.11,29.12){1/2m}
\rput(104.12,116.13){1/q=1/2+1/p'}
\rput(78.61,161.97){1/q=(m+1)/p'}
\rput(120.24,77.8){1/q=3/p'}
\psline[linewidth=0.6](24.19,37.07)(50.79,94.79)
\psline[linewidth=0.55,linestyle=dashed,dash=1 1](50.79,94.79)(71.78,142.14)
\newrgbcolor{userFillColour}{1 0 0.8}
\psline[linewidth=0.4,fillcolor=userFillColour,fillstyle=solid,arrowsize=1.3 2]{->}(24.19,36.99)(24.19,150.4)
\newrgbcolor{userFillColour}{1 0 0.8}
\psline[linewidth=0.4,fillcolor=userFillColour,fillstyle=solid](24.19,89.65)(182.21,121.24)
\newrgbcolor{userFillColour}{1 0 0.8}
\psline[linewidth=0.4,fillcolor=userFillColour,fillstyle=solid](155.87,115.97)(24.19,36.99)
\rput(15.74,141.68){1}
\rput(79.47,7.37){\textbf{Figure \thesection.\theabb:} Range of p' and q}
\end{pspicture}

The proof will be the main work in this paper and be done in the next chapters. First of all, we will derive the restriction theorem.\\
As an immediate consequence of Theorem \ref{superlemma}, we get as corollary

\begin{kornr}
Let $p'>m+1$, $\frac{3}{p'}<\tfrac{1}{q}$. Then 
\addtocounter{equation}{-1}
\begin{align}
		\|\sum_{\alpha} a_\alpha \hat\phi_\alpha\|_{p'}	\lesssim&  
		\delta^{\durch{q}}\|\sum_{\alpha=(k,j,n)} 
							a_\alpha (\zweik\gamma)^{\durch{q}-\frac{m+1}{p'}} \phi_\alpha\|_{q',1}.
\end{align}
\end{kornr}


\begin{bew}
Actually, we will just use that Theorem \ref{superlemma} is valid for every $p'>m+1$ and for some range $\frac{3}{p'}<\tfrac{1}{q}<\frac{3}{p'}+\e_p$. According to the assumptions, we have $p'>m+1\geq4$.
Hence $\frac{3}{p'}=\frac{2}{p'}+\durch{p'}<\durch{2}+\durch{p'}$.
Therefore we find $r$ with $\frac{3}{p'}<\durch{r}<\frac{1}{2}+\frac{1}{p'}$ (cf. Figure \thesection.\theabb), i.e. satisfying the requirements of Theorem \ref{superlemma}.
So it is sufficient to show that if \eqref{arclength} holds for some $(p,r)$, then it also holds for $(p,q)$ with $\durch{r}<\durch{q}$. 
Under this condition, there exists $1\leq s<\infty$ such that $\durch{r'}=\durch{q'}+\durch{s}$.
This means
\begin{align}\label{rqs}
	\durch{r}-\durch{q}=\durch{q'}-\durch{r'}=-\durch{s}.
\end{align}
For $(x,y,z)$ in the support of $\phi_{kjn}$, it follows $x\approx\zweik\gamma$, so we introduce 
$g(x,y,z):=x^{-\durch{s}}$. We now may apply Theorem \ref{superlemma} and use Hölder's inequality for Lorentz spaces (see Lemma \ref{lohoe}):
\begin{align}\label{satzanwenden}
		\|\sum_{\alpha} a_\alpha \hat\phi_\alpha\|_{p'}	\lesssim&  
		\delta^{\durch{r}} \|\sum_{\alpha} 
							a_\alpha (\zweik\gamma)^{\durch{r}-\frac{m+1}{p'}} \phi_\alpha\|_{r'} \nonumber\\
		=& \delta^{\durch{r}} \|\sum_{\alpha} 
					a_\alpha (\zweik\gamma)^{\durch{q}-\frac{m+1}{p'}-\durch{s}} \phi_\alpha\|_{r'}\nonumber\\
		\lesssim& \delta^{\durch{r}} \|\sum_{\alpha} a_\alpha (\zweik\gamma)^{\durch{q}-\frac{m+1}{p'}}
																 \phi_\alpha g \|_{L^{r'}(\Gamma^\delta)}	\nonumber \\
		\leq& \delta^{\durch{r}} \|\sum_{\alpha} 	
		a_\alpha (\zweik\gamma)^{\durch{q}-\frac{m+1}{p'}} \phi_\alpha \|_{L^{q',r'}(\Gamma^\delta)}
						 \|g\|_{L^{s,\infty}(\Gamma^\delta)}.				
\end{align}
A short computation yields
\begin{align*}
	|\{(x,y,z)\in\Gamma^\delta: x^{-\durch{s}} >\lambda \}| =&
	\int_0^{\lambda^{-s}} \int_1^2 \int_{z\Phi(\tfrac{x}{z})}^{z\Phi(\tfrac{x}{z})+z\delta} 
																										\dd y\dd z\dd x \\
	=& \delta \lambda^{-s} \int_1^2 z\dd z \approx \delta \lambda^{-s},															
\end{align*}
i.e.
\begin{align}\label{20OCT11}
	\|g\|_{L^{s,\infty}(\Gamma^\delta)}=\sup_{\lambda>0} \lambda |\{(x,y,z)\in\Gamma^\delta: x^{-\durch{s}} >\lambda \}|^{\durch{s}} \approx \delta^{\durch{s}}.
\end{align}
Combined with \eqref{satzanwenden}, we end up with
\begin{align*}
		\|\sum_\alpha a_{\alpha} \hat\phi_\alpha\|_{p'}	\lesssim&  
		\delta^{\durch{r}+\durch{s}} \|\sum_{\alpha} 	a_\alpha (\zweik\gamma)^{\durch{q}-\frac{m+1}{p'}}
																	 \phi_\alpha \|_{q',r'}\\
		\leq& \delta^{\durch{q}} \|\sum_{\alpha} 	a_\alpha (\zweik\gamma)^{\durch{q}-\frac{m+1}{p'}}
																					\phi_\alpha \|_{q',1}\ ,
\end{align*}
since \eqref{rqs} means $\durch{r}+\durch{s}=\durch{q}$.
\end{bew}

For $\durch{q}=\frac{m+1}{p'}$ as considered in Theorem \ref{thm1}, \eqref{arclength} reads
\begin{align}
		\|\sum_\alpha a_{\alpha} \hat\phi_\alpha\|_{p'}	\lesssim&  
		\delta^{\durch{q}} \|\sum_{\alpha} 	a_\alpha \phi_\alpha \|_{(q',1)}
\end{align}
for all sequences $a_\alpha$ and every $\delta>0$. 
This is nothing but a "desingularised" adjoint restriction estimate with the singular surface measure replaced by a $\delta$-thickening of the surface, for certain discrete functions, and it implies the desired restriction theorem. To be more precise, for 
$\durch{q}=\frac{m+1}{p'}$ and $q<\infty$ or $p'> m+1$ respectively, the weak $L^p(\R^3)$-$L^{q,\infty}(\Gamma)$ estimate
\begin{align*}
		\|\hat f|_\Gamma \|_{L^{q,\infty}(\Gamma)}	\lesssim& \|f \|_p\qquad\forall f\in\s(\R^3)
\end{align*}
holds. Eventually we finish the proof by using the generalised Marcinkiewicz interpolation theorem (for instance, see Theorem 1.4.19 in [G]).


\section{Estimation of the overlap}\label{kap:ueberlapp}

\subsection{Straightforward results}
It will become essential for us to understand the overlap of the sets  $\Gamma_{kjn}^\delta+\Gamma_{lip}^\delta$, $k,l=\krange$, $j\in I_k$, $i\in I_l$, $n,p=\zrange$, i.e. to examine the maximal number of them containing a single point. However, this number will not be bounded by a absolute constant (which would be optimal, but is not true). Instead will collect some overlap from the $z-$direction, which however is natural and manageable.

\begin{lemnr}\label{ueberlappz}
For every $\xi\in\R^3$ we have
\begin{align}
\#\{(n,p)|\exists k,l,j,i:\xi\in\Gammaplus\}\leq 3 \beta\inv.
\end{align}
\end{lemnr}~\\
\begin{bew}
If $\xi\in\Gammaplus$, there exist $\textbf{x}_1\in\Gamma_{kjn}^\delta$ and $\textbf{x}_2\in\Gamma_{lip}^\delta$ with $\xi=\textbf{x}_1+\textbf{x}_2$. Let $z_1,z_2$ denote the last component of $\textbf{x}_1$ and $\textbf{x}_2$ respectively, which means
$n\beta\leq z_1\leq(n+1)\beta$ and $p\beta\leq z_2\leq(p+1)\beta$. The sum of both inequalities leads to $(n+p)\beta\leq z\leq (n+p+2)\beta$. This implies
$\frac{z}{\beta}-2\leq n+p \leq \frac{z}{\beta}$ and hence
\begin{align*} n\in -p+[\tfrac{z}{\beta}-2,\tfrac{z}{\beta}]\cap\N.\end{align*}
We conclude
\begin{align*}
	\#\{(n,p)|\exists k,l,j,i:\xi\in\Gammaplus\}
	=& 		\sum_{p=\durch{\beta}}^{\frac{2}{\beta}-1} \#\{n|\exists k,l,j,i:\xi\in\Gammaplus\} \\
	\leq& \sum_{p=\durch{\beta}}^{\frac{2}{\beta}-1}
									 \#\left([\tfrac{z}{\beta}-2,\tfrac{z}{\beta}]\cap\N\right) \\
	\leq& \sum_{p=\durch{\beta}}^{\frac{2}{\beta}-1} 3 = 3 \beta\inv,
\end{align*}
completing the proof.
\end{bew}
~\\

Exploiting the dyadic structure of the decomposition, we obtain a further simple result:\\
\begin{lemnr}\label{ueberlappk}
For $k=\krange$ let $V_k(n,p)=\bigcup\{\Gamma_{kjn}^\delta+\Gamma_{lip}^\delta|l\leq k,\,j\in I_k,\,i\in I_l\}$. Then for all $\xi\in\R^3$, $n,p=\zrange$ holds
\begin{align*}
	\#\{k|\xi\in V_k(n,p)\}\leq 3.
\end{align*}
\end{lemnr}
\begin{bew}
The procedure is comparable to the previous proof, though we now concentrate on the $x$-component. If $\xi\in V_k(n,p)=\bigcup\{\Gamma_{kjn}^\delta+\Gamma_{lip}^\delta|l\leq k,\,j\in I_k,\,i\in I_l\}$ and $x$ denotes the first component of $\xi$, there exist $l,j,i$ such that
$$\zweik\gamma+j\gamma_k+\zweil\gamma+i\gamma_l \leq x\leq \zweik\gamma+(j+1)\gamma_k+\zweil\gamma+(i+1)\gamma_l.$$
Since we are just interested in the dyadic size, we estimate quite roughly
$$\zweik\gamma \leq x\leq \zweik\gamma+2^{k\mhalb}2^{k(1-\mhalb)}\gamma +\zweil\gamma+2^{l\mhalb}2^{l(1-\mhalb)}\gamma= (2^{k+1}\minus 1)\gamma+(2^{l+1}\minus 1)\gamma.$$
and since we assumed $l\leq k$, this reduces to
$$\zweik\gamma \leq x\leq (2^{k+1}\minus 1)\gamma+(2^{l+1}\minus 1)\gamma \leq (2^{k+2}\minus 1)\gamma.$$
Hence
$$2^k\leq\frac{x}{\gamma}+1\leq2^{k+2}$$
and therefore
$$\log_2\left(\frac{x}{\gamma}+1\right)-2\leq k\leq \log_2\left(\frac{x}{\gamma}+1\right).$$
Our claim is an immediate consequence of this inequality.
\end{bew}
~\\

The handling of the overlap concerning the remaining parameters is more complicated. Here we need to involve the $y-$coordinate. Therefore, we need as a start a new coordinate system.\\

\subsection{Further results}


\begin{lemnr} \label{ueberlapprest}
There exists an absolute constant $C>0$ (not depending on $\delta$) such that
for all $\xi\in\R^3$, $k=\krange$, $n,p=\zrange$
\begin{align}
	\#\{(j,l,i)|l\leq k,\xi\in\Gammaplus\}\leq C.
\end{align}
\end{lemnr}~\\
\begin{bew}
We are allowed to restrict ourself to the case $l\ll k$ (more precise: $l<k-2$), since in the case $l\approx k$ the second derivative of $\Phi$ and therefore the Gaussian curvature is comparable on the regions $x\approx 2^k\gamma$ and $x\approx2^l\gamma$. After rescaling, we are back in the classical case with (almost) constant curvature.
\\
The procedure differs a bit from the previous lemmas:
We claim that the number of triples $(j,l,i)$, $l\leq k$, with the property $\xi\in\Gamma_{kjn}+\Gamma_{lip}^\delta$ is bounded by a fixed constant. So let
\begin{eqnarray}\label{siemeinen1s26}	\xi\in(\Gamma_{kjn}+\Gamma_{lip}^\delta)\cap(\Gamma_{kj'n}+\Gamma_{l'i'p}^\delta),
\end{eqnarray}
$l,l'\leq k$.
Without loss of generality, we may assume $l'\leq l$ and in the case $l'=l$ furthermore $i'\leq i$ by interchanging the parameters if necessary.

\begin{stepnr}\label{step1}
If $(l',i')\neq(l,i)$ then $j\leq j'$.
\end{stepnr}
The assumption is not necessary, since in the case $(l',i')=(l,i)$ we may ensure $j\leq j'$ by interchanging the parameters with and without primes.

\underline{Case 1: $l'\neq l$}
\begin{quote}
In this case is
 $l'<l$, i.e. $l'+1\leq l$, and therefore
\begin{align*}
	& \zweik\gamma+j\gamma_k+\zweil\gamma \\
	\leq& \zweik\gamma+j\gamma_k+\zweil\gamma +i\gamma_l \leq x  \\
	<& \zweik\gamma+(j'+1)\gamma_k+(2^{l'}\minus1)\gamma +(i'+1)\gamma_{l'} \\
	\leq& \zweik\gamma+(j'+1)\gamma_k+(2^{l'+1}\minus1)\gamma								\\
	\leq& \zweik\gamma+(j'+1)\gamma_k+(2^{l}\minus1)\gamma .
\end{align*}
This already implies $j< j'+1$, i.e. $j\leq j'$.
\end{quote}
\underline{Case 2: $l'= l$}
\begin{quote}
In this case is $i'<i$, i.e. $i'+1\leq i$ and therefore
\begin{align*}
			& \zweik\gamma+j\gamma_k+\zweil\gamma +i\gamma_l \leq x  \\
	<& \zweik\gamma+(j'+1)\gamma_k+\zweil\gamma +(i'+1)\gamma_{l} \\
	\leq& \zweik\gamma+(j'+1)\gamma_k+\zweil\gamma+i\gamma_{l}.							
\end{align*}
This again implies $j< j'+1$, i.e. $j\leq j'$.
\end{quote}

In the next step, it is useful to examine the projections in $x$-$y-$space. Thus we introduce the new curves $\Phi_n(x)=n\beta\Phi\left(\frac{x}{n\beta}\right)$ (analogue $\Phi_p$). 
Using $\beta\leq\delta$, it is an easy task to verify
\begin{align}\label{phin}
	\Gamma_{kjn}^\delta=&\{(x,y,z)\in\R^3|\, n\beta\leq z\leq(n+1)\beta,~
	\Phi\left(\tfrac{x}{z}\right)\leq\tfrac{y}{z}\leq\Phi\left(\tfrac{x}{z}\right)+\delta,~
  x_{kj}\leq x\leq x_{k,j+1}\}  \\
  \subset \tilde\Gamma_{kjn}^{\delta}=& \{(x,y,z)\in\R^3|\, n\beta\leq z\leq(n+1)\beta,~
	|y-\Phi_n(x)|\leq 10m\delta,~
  x_{kj}\leq x\leq x_{k,j+1}\},
\end{align}
and likewise for $\Gamma_{lip}^\delta$.
Now define the projection $P(x,y,z)=y-x\Phi_n'(x_k)$ on the normal to the graph of $\Phi_n$ at point $x_k$ in $x$-$y-$space.\\

\stepcounter{abb}
\def\JPicScale{0.5}
\psset{unit=\JPicScale mm}
\psset{linewidth=0.3,dotsep=1,hatchwidth=0.3,hatchsep=1.5,shadowsize=1,dimen=middle}
\psset{dotsize=0.7 2.5,dotscale=1 1,fillcolor=black}
\psset{arrowsize=1 2,arrowlength=1,arrowinset=0.25,tbarsize=0.7 5,bracketlength=0.15,rbracketlength=0.15}
\hspace{20mm}\begin{pspicture}(0,0)(157.5,90)
\newrgbcolor{userFillColour}{0.87 0.85 0.85}
\pspolygon[linestyle=none,fillcolor=userFillColour,fillstyle=solid](73.82,17.86)
(73.82,12.47)
(113.16,46.15)
(106.71,53.65)
(59.87,13.65)(73.82,17.86)
\newrgbcolor{userFillColour}{0.6 0.6 0.6}
\pspolygon[linestyle=none,fillcolor=userFillColour,fillstyle=solid](60.13,8.52)(73.82,12.47)(73.82,17.73)(60.13,13.78)
\rput(60.52,-4.74){$x_{li}$}
\psline{<-}(150,0)(0,0)
\psline{->}(0,0)(0,90)
\rput(115,3.75){$x_k$}
\rput(123.75,-5){$x_{kj}$}
\rput(141.25,-5){$x_{k,j+1}$}
\rput(17.63,-5){$x_{l'i'}$}
\rput(143.75,42.5){$\tilde\Gamma_{kjn}^\delta$}
\rput(16.58,13.16){$\tilde\Gamma_{l'i'p}^\delta$}
\rput(62.37,21.71){$\tilde\Gamma_{lip}^\delta$}
\psbezier(150.13,79.87)(124.34,0)(0,0.13)(0,0.13)
\psbezier(150,85)(124.21,5.13)(-0.13,5.26)(-0.13,5.26)
\psline(20.13,1.38)(20.13,6.58)
\psline(35,3.16)(35,8.42)
\psline(60.13,8.55)(60.13,13.82)
\psline(73.82,12.63)(73.82,17.76)
\psline{->}(120.53,37.37)(94.34,68.55)
\psline(125.92,42.11)(125.92,47.37)
\psline(135.26,52.24)(135.26,57.5)
\psline[linestyle=dotted](125.92,41.97)(125.92,0.13)
\psline[linestyle=dotted](135.26,52.76)(135.26,0)
\psline{[-]}(125.92,0)(135.26,0)
\psline{[-]}(60,0)(73.95,0)
\psline{[-]}(20,0)(35.13,0)
\psline(150.13,79.87)(150.26,85.26)
\rput(153.03,82.5){$\delta$}
\psline[linestyle=dotted](73.95,12.76)(113.16,46.05)
\psline[linestyle=dotted](60.13,13.95)(106.97,53.95)
\psline{[-]}(113.16,46.05)(106.71,53.68)
\rput(115,60){$P(\tilde\Gamma_{lip}^\delta)$}
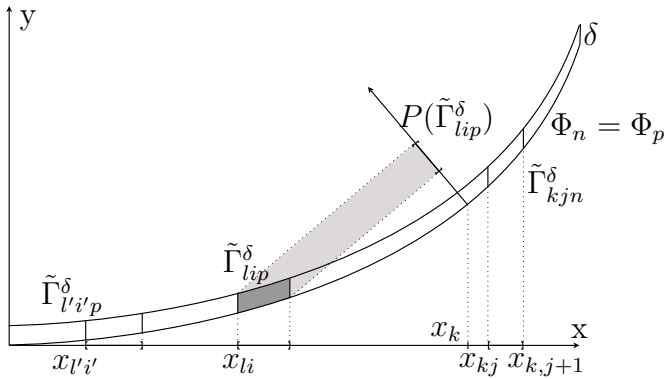
\rput(70,-20){\textbf{Figure \thesection.\theabb:} Projection on the normal vector in the case $n=p$}
\rput(150,3.75){x}
\rput(5,86.25){y}
\psline[linestyle=dotted](120.66,37.34)(120.66,-0.03)
\psline(120.66,-0.95)(120.66,0.89)
\psline[linestyle=dotted](60.13,-0.03)(60.13,8.52)
\psline[linestyle=dotted](73.82,-0.03)(73.82,12.73)
\psline[linestyle=dotted](34.87,-0.16)(35,3.26)
\rput(157.5,57.5){$\Phi_n=\Phi_p$}
\end{pspicture}
~\\
\vspace{5mm}

Additionally, we introduce
\begin{align}\label{xsidef}
	\xi_{kjn}=(x_{kj},\Phi_n(x_{kj}),n\beta)\in\Gamma_{kjn}^\delta\qquad \xi_{lip}=(x_{li},\Phi_p(x_{li}),p\beta)\in\Gamma_{lip}^\delta
\end{align}
 and analogue $\xi_{kj'n}$, $\xi_{l'i'p}$.

\begin{stepnr}\label{step2}~
\begin{enumerate} \renewcommand{\labelenumi}{(\roman{enumi})}
	\item $P(\xi_{l'i'p})-P(\xi_{lip}) \approx \delta 2^{k(m-1)-l'(\mhalb-1)} \Delta_{ll'ii'}$
	\item $P(\xi_{kj'n})-P(\xi_{kjn}) \gtrsim\delta j(j'-j),\qquad\text{ if }j\leq j'+1$   
\end{enumerate} 
with $\Delta_{ll'ii'}=2^{l'(\mhalb-1)}\dfrac{x_{li}-x_{l'i'}}{\gamma}=2^{l'(\mhalb-1)}
			\left(2^l+i2^{l(1-\mhalb)}-2^{l'}-i'2^{l'(1-\mhalb)}\right)$.
\end{stepnr}
The term $\Delta_{ll'ii'}$ looks somehow artificial. However, we will discover that this quantity is the crucial one, expressing (in some sense) the distance between $\Gamma_{lip}^\delta$ and $\Gamma_{l'i'p}^\delta$.\\
\underline{Part (i):}
The mean value theorem provides the existence of an $\tilde x\in (x_{l'i'},x_{li})$ with 
\begin{align}\label{step2a}\Phi_p(x_{l'i'})-\Phi_p(x_{li})=\Phi_p'(\tilde x)(x_{l'i'}-x_{li}).
\end{align}
Thus especially
\begin{align}\nonumber
	\tilde x\leq x_{li}\leq x_{l+1}\underset{l<k-2}{\leq} x_{k-2}
\end{align}
and hence
\begin{align*}
	\frac{\tilde x}{p\beta} \leq \frac{n\beta}{p\beta}\frac{x_{k-2}}{n\beta} 
	\leq 	\frac{2x_{k-2}}{n\beta} 
	= \frac{2(2^{k-2}\minus1)\gamma}{n\beta} 
	= \halb \frac{(2^{k}\minus4)\gamma}{n\beta}
	\leq \halb \frac{(2^{k}\minus1)\gamma}{n\beta}
	= \halb \frac{x_k}{n\beta}.
\end{align*}
Using monotony and convexity of $\Phi'$, as well as $\Phi'(0)=0$, we obtain
\begin{align}\label{step2b} 
		\Phi'\left(\frac{\tilde x}{p\beta}\right) 
		\leq \Phi'\left(\halb \frac{x_k}{n\beta}\right) 
		\leq \halb \Phi'\left(\frac{x_k}{n\beta}\right).		
\end{align}
Hence
\begin{eqnarray}
	\Phi_n'(x_k)-\Phi_p'(\tilde x)
	&=& \Phi'\left(\frac{x_k}{n\beta}\right)-\Phi'\left(\frac{\tilde x}{p\beta}\right) \nonumber\\
	&\underset{\eqref{step2b}}{\approx}& \Phi'\left(\frac{x_k}{n\beta}\right) \nonumber\\
	&\approx& \frac{x_k^{m-1}}{(n\beta)^{m-1}} \nonumber\\
	&\approx& \label{step2c} \gamma^{m-1}2^{k(m-1)}.
\end{eqnarray}
Thus we conclude
\begin{eqnarray*}
	P(\xi_{l'i'p})-P(\xi_{lip}) &=& \Phi_p(x_{l'i'})-\Phi_p(x_{li})-(x_{l'i'}-x_{li})\Phi_n'(x_k)\\
	&\underset{\eqref{step2a}}{=}& (\Phi_n'(x_k)-\Phi_p'(\tilde x))(x_{li}-x_{l'i'}) \\
	&=& (\Phi_n'(x_k)-\Phi_p'(\tilde x))2^{-l'(\mhalb-1)}\Delta_{ll'ii'}\gamma	\\
	&\underset{\eqref{step2c}}{\approx}& \gamma^m 2^{k(m-1)} 2^{-l'(\mhalb-1)}\Delta_{ll'ii'} \\
	&=&\delta 2^{k(m-1)-l'(\mhalb-1)}\Delta_{ll'ii'}
\end{eqnarray*}

\underline{Part (ii):} \\
Since the projection $P$ depends on $k$, we are not able to use (i). However, we proceed in a similar manner: like in part (i), we obtain a $\tilde x_k$ between $x_{kj}$ and $x_{kj'}$ with
\begin{align*}
	\Phi_n(x_{kj'})-\Phi_n(x_{kj})=\Phi_n'(\tilde x_k)(x_{kj'}-x_{kj}).
\end{align*}
Furthermore, we obtain a $\tilde{\tilde {x}}_k$ between $\tilde x_k$ and $x_k$ with
\begin{align}\label{step2l}
	P(\xi_{kj'n})-P(\xi_{kjn}) =& \Phi_n(x_{kj'})-\Phi_n(x_{kj})-(x_{kj'}-x_{kj})\Phi_n'(x_k)
																													 \nonumber\\
	{=}& (\Phi_n'(\tilde x_k)-\Phi_n'(x_k))(x_{kj'}-x_{kj}) \\
	=& \Phi_n''(\tilde{\tilde {x}}_k)(\tilde x_k-x_k)(x_{kj'}-x_{kj}).\nonumber
\end{align}
Especially $\tilde{\tilde {x}}_k\in(x_k,\tilde x_k)\subset(x_k,x_{k+1})$, i.e.
\begin{align}\label{step2m}
	\tilde{\tilde {x}}_k \approx 2^k\gamma.
\end{align}

In the case $j'\geq j$ we use
\begin{align}\label{step2na}
	\tilde x_k > \min\{x_{kj},x_{kj'}\}=x_{kj}= \zweik\gamma+j2^{k(1-\mhalb)}\gamma = x_k+j2^{k(1-\mhalb)}\gamma.
\end{align}
This gives
\begin{eqnarray*}
	P(\xi_{kj'n})-P(\xi_{kjn}) 
	&\underset{\eqref{step2l}}{=}& \Phi_n''(\tilde{\tilde {x}}_k)(\tilde x_k-x_k)(x_{kj'}-x_{kj}) \\
	&\underset{\eqref{step2m}}{\approx}& (2^k\gamma)^{m-2}(\tilde x_k-x_k)(x_{kj'}-x_{kj}) \\
	&\underset{\eqref{step2na}}{\geq}& (2^k\gamma)^{m-2} ~ j2^{k(1-\mhalb)}\gamma ~ 	
											(j'-j)2^{k(1-\mhalb)}\gamma	\\
	&=& \gamma^m\, j(j'-j) \\
	&=& \delta j(j'-j).
\end{eqnarray*}
In the other case $j'<j$ the assumption ensures $j=j'+1$, i.e. $j'-j=-1$. Now we use
\begin{align}\label{step2nb}
	\tilde x_k < \max\{x_{kj},x_{kj'}\} = x_{kj}= \zweik\gamma+j2^{k(1-\mhalb)}\gamma = x_k+j2^{k(1-\mhalb)}\gamma.
\end{align}
In this case we obtain
\begin{eqnarray*}
	P(\xi_{kj'n})-P(\xi_{kjn}) 
	&\underset{\eqref{step2l}}{=}& \Phi_n''(\tilde{\tilde {x}}_k)(\tilde x_k-x_k)(x_{kj'}-x_{kj}) \\
	&\underset{\eqref{step2m}}{\approx}& (2^k\gamma)^{m-2}(\tilde x_k-x_k)(-2^{k(1-\mhalb)}\gamma) \\
	&\underset{\eqref{step2nb}}{\geq}& -(2^k\gamma)^{m-2} ~ j2^{k(1-\mhalb)}\gamma ~ 	
											2^{k(1-\mhalb)}\gamma	\\
	&=& -\gamma^m\, j \\
	&=& -\delta j \ =\ \delta j(j'-j),
\end{eqnarray*}
completing Step \ref{step2}.

The next task will be to estimate the size of the pieces $\Gamma_{kjn}^\delta$ of our decomposition,
with respect to the projection $P$. The size of a set $U$ with respect to $P$ is measured by
$\diam_P(U)=\sup\{|P(u)-P(v)|:u,v\in U\}$.

\begin{stepnr}\label{step3}~
\begin{enumerate}\renewcommand{\labelenumi}{(\roman{enumi})}
	\item $\diam_P(\Gamma_{lip}^\delta) \lesssim \delta 2^{k(m-1)-l(\mhalb-1)}$
	\item $\diam_P(\Gamma_{kjn}^\delta) \lesssim \delta 2^{k\mhalb}$
\end{enumerate}
\end{stepnr}

\underline{Part (i):}
\begin{align*}
	\diam_P(\Gamma_{lip}^\delta) \underset{}{=}&
				P(x_{li},\Phi_p(x_{li})+\landau(\delta),p\beta)-P(x_{l,i+1},\Phi_p(x_{l,i+1}),p\beta) \\
	=& P(x_{li},\Phi_p(x_{li}),p\beta)-P(x_{l,i+1},\Phi_p(x_{l,i+1}),p\beta)+\landau(\delta) \\
	=& P(\xi_{lip})-P(\xi_{l,i+1,p})+\landau(\delta)
\end{align*}
We apply Step \ref{step2}(i), replacing $l'$ by $l$, $i$ by $i+1$ and $i'$ by $i$.
Then\\ $\Delta_{ll,i+1,i}=2^{l(\mhalb-1)}	\left((i+1)2^{l(1-\mhalb)}-i2^{l(1-\mhalb)}\right)=1$ and
\begin{align*}
	\diam_P(\Gamma_{lip}^\delta)
	=& P(\xi_{lip})-P(\xi_{l,i+1,p})+\landau(\delta) \\
	\approx& \delta 2^{k(m-1)-l(\mhalb-1)} \Delta_{ll,i+1,i} +\landau(\delta) 	\\
	=& \delta 2^{k(m-1)-l(\mhalb-1)} +\landau(\delta) \\
	\overset{l\leq k}{\lesssim}& \delta 2^{k(m-1)-l(\mhalb-1)}.
\end{align*}

\underline{Part (ii):}
\begin{align*}
	\diam_P(\Gamma_{kjn}^\delta) \underset{}{=}&
				P(x_{k,j+1},\Phi_n(x_{k,j+1})+\landau(\delta),n\beta)-P(x_{kj},\Phi_n(x_{kj}),n\beta) \\
	=& P(x_{k,j+1},\Phi_n(x_{k,j+1}),n\beta)-P(x_{kj},\Phi_n(x_{kj}),n\beta)+\landau(\delta) \\
	=& P(\xi_{k,j+1,n})-P(\xi_{kjn})+\landau(\delta)
\end{align*}
Step \ref{step2}(ii) implies
$$P(\xi_{kj'n})-P(\xi_{kjn}) \gtrsim  \delta j(j'-j),\qquad\text{ if }j\leq j'+1$$
or, equivalently
$$P(\xi_{kjn})-P(\xi_{kj'n}) \lesssim \delta j(j-j'),\qquad\text{ if }j\leq j'+1,$$
so especially
$$P(\xi_{k,j+1n})-P(\xi_{kjn}) \lesssim (j+1)\delta\leq \delta 2^{k\mhalb}.$$
This leads to
\begin{align*}
	\diam_P(\Gamma_{kjn}^\delta) 
	=& P(\xi_{k,j+1,n})-P(\xi_{kjn})+\landau(\delta) \\
	\lesssim& \delta 2^{k\mhalb} +\landau(\delta) 	\\
	\lesssim& \delta 2^{k\mhalb}.
\end{align*}~

A further step will analyse the quantity $\Delta_{ll'ii'}$. Here we exploit \eqref{siemeinen1s26}: $\Gamma_{kj'n}^\delta+\Gamma_{l'i'p}^\delta\neq$\O\ .

\begin{stepnr}\label{step4}
	$\Delta_{ll'ii'} \lesssim 1$
\end{stepnr}
We apply Lemma \ref{überlapplemma}. Since
$$\xi_\alpha\in\Gamma_\alpha^\delta,\qquad \alpha=(k,j,n),(k,j',n),(l,i,p),(l',i',p)$$
the lemma yields
\begin{eqnarray}\label{blablub}
	& &P(\xi_{kj'n})-P(\xi_{kjn}) + P(\xi_{l'i'p})-P(\xi_{lip}) \nonumber\\
	& \leq&\diam_P(\Gamma_{kjn}^\delta)+\diam_P(\Gamma_{kj'n}^\delta)
	+\diam_P(\Gamma_{lip}^\delta)+\diam_P(\Gamma_{l'i'p}^\delta).
\end{eqnarray}
In Step \ref{step1}, we established $j\leq j'$ (at least if $(l,i)\neq(l',i')$, but in the case $(l,i)=(l',i')$ we may also assume this). The application of Step \ref{step2} and Step \ref{step3} results in
\begin{align}\label{step3a}
	\delta 2^{k(m-1)-l'(\mhalb-1)} \Delta_{ll'ii'} \leq&
	\delta 2^{k(m-1)-l'(\mhalb-1)} \Delta_{ll'ii'} + \delta j(j'-j)  \nonumber\\
	\overset{\eqref{blablub}}{\lesssim}& \delta 2^{k\mhalb} + \delta 2^{k(m-1)-l(\mhalb-1)}+\delta 2^{k(m-1)-l'(\mhalb-1)}
\end{align}
Now we use $l'\leq l$ and $l'<k$, since the last one implies
$k(m-1)-l'(\mhalb-1)=k\mhalb+(k-l')(\mhalb-1)\geq k\mhalb$. Thus \eqref{step3a} transforms into
\begin{align*}
	2^{k(m-1)-l'(\mhalb-1)} \Delta_{ll'ii'}  
	\lesssim  2^{k(m-1)-l(\mhalb-1)}+2^{k(m-1)-l'(\mhalb-1)}+ 2^{k\mhalb} \lesssim 2^{k(m-1)-l'(\mhalb-1)},
\end{align*}
i.e. $\Delta_{ll'ii'} \lesssim 1$.

The former results will be merged into the following step, which states that there are not "too many" $(l',i')$ appropriate to a given $(l,i)$.
\begin{stepnr}\label{step5}
	One of these three alternatives holds:
\begin{enumerate}
	\item $l'=l$: Then $|i-i'|\lesssim1$.
	\item $l'=l-1$: Then $i\lesssim1$ and $2^{l'\mhalb}-i'\lesssim1$.
	\item $l'<l-1$: Then $l',l\lesssim1$ and especially $i',i\lesssim1$.
\end{enumerate}
Especially $l'\approx l$ is necessary.
\end{stepnr}
Let me remind you that we assumed $l'\leq l$ and even $i'\leq i$ for $l'=l$. The proof is divided into three cases:\\
\underline{Case 1: $l'=l$}\\
According to Step \ref{step4}, we know
$$1\gtrsim\Delta_{llii'} 
	= 2^{l(\mhalb-1)}	\left(2^l+i2^{l(1-\mhalb)}-2^{l}-i'2^{l(1-\mhalb)}\right) 
	= i-i' \geq 0.$$

\underline{Case 2: $l'=l-1$}\\
Step \ref{step4} now reads
\begin{align*}
		1\gtrsim\Delta_{ll'ii'} 
		= 2^{l'(\mhalb-1)}	\left(2^{l'+1} +i2^{(l'+1)(1-\mhalb)} -2^{l'} -i'2^{l'(1-\mhalb)}\right) 
		= \underbrace{2^{l'\mhalb}-i'}_{\geq0} + \underbrace{i2^{1-\mhalb}}_{\geq0},
\end{align*}
hence
\begin{align*}
	1 \gtrsim i2^{1-\mhalb} \approx i \qquad\text{and}\qquad 1 \gtrsim  2^{l'\mhalb} -i'.
\end{align*}

\underline{Case 3: $l>l'+1$}\\
In this case, Step \ref{step4} yields
\begin{align*}
	1\gtrsim\Delta_{ll'ii'} 
		=& 2^{l'(\mhalb-1)}	\left(2^l+i2^{l(1-\mhalb)}-2^{l'}-i'2^{l'(1-\mhalb)}\right) \\
		\gtrsim& 2^{l'(\mhalb-1)}	\left(2^l -2^{l'+1}\right) 
		\gtrsim 2^{l'(\mhalb-1)} 2^l \geq	2^l,
\end{align*}
which already implies $l'\leq l\lesssim 1$ and thus 
$i\leq 2^{l\mhalb} \lesssim 1$, as well as $i'\lesssim 1$.\\

\stepcounter{abb}
\psset{unit=\JPicScale mm}
\psset{linewidth=0.3,dotsep=1,hatchwidth=0.3,hatchsep=1.5,shadowsize=1,dimen=middle}
\psset{dotsize=0.7 2.5,dotscale=1 1,fillcolor=black}
\psset{arrowsize=1 2,arrowlength=1,arrowinset=0.25,tbarsize=0.7 5,bracketlength=0.15,rbracketlength=0.15}
\hspace{20mm}\begin{pspicture}(0,0)(158.03,85)
\newrgbcolor{userFillColour}{0.89 0.88 0.88}
\pspolygon[linestyle=none,fillcolor=userFillColour,fillstyle=solid](130.79,37.24)
(140.13,47.24)
(140.13,52.5)
(4.87,52.5)
(4.87,37.24)(130.79,37.24)
\newrgbcolor{userFillColour}{0.6 0.6 0.6}
\pspolygon[linestyle=none,fillcolor=userFillColour,fillstyle=solid](130.92,37.11)
(140.26,46.97)
(140.26,52.37)
(130.92,42.24)(130.92,37.11)
\rput(43.75,-10){$x_{li}$}
\psline{<-}(155,-5)(5,-5)
\psline{->}(5,-5)(5,85)
\rput(130,-10){$x_{kj}$}
\rput(147.5,-10){$x_{k,j+1}$}
\rput(21.25,-10){$x_{l'i'}$}
\rput(147.5,38.75){$\tilde\Gamma_{kjn}^\delta$}
\rput(21.25,6.25){$\tilde\Gamma_{l'i'p}^\delta$}
\rput(43.75,10){$\tilde\Gamma_{lip}^\delta$}
\psbezier(155.13,74.87)(129.34,-5)(5,-4.87)(5,-4.87)
\psbezier(155,80)(129.21,0.13)(4.87,0.26)(4.87,0.26)
\psline(25.13,-3.62)(25.13,1.58)
\psline(40,-1.84)(40,3.42)
\psline(53.75,1.12)(53.75,6.25)
\psline(130.92,37.11)(130.92,42.37)
\psline(140.26,47.24)(140.26,52.5)
\psline[linestyle=dotted](130.92,36.97)(130.92,-4.87)
\psline[linestyle=dotted](140.26,47.76)(140.26,-5)
\psline{[-]}(130.92,-5)(140.26,-5)
\psline{[-]}(40,-5)(53.95,-5)
\psline{[-]}(25,-5)(40.13,-5)
\psline(155.13,74.87)(155.26,80.26)
\rput(158.03,77.5){$\delta$}
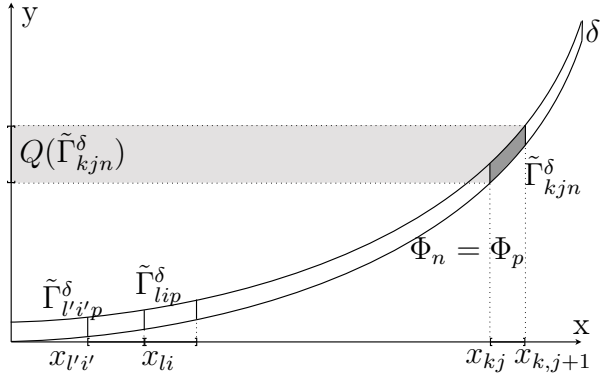
\rput(70,-25){\textbf{Figure \thesection.\theabb:} Projection on the y-axis in the case $n=p$}
\rput(155,-1.25){x}
\rput(10,81.25){y}
\psline[linestyle=dotted](39.87,-5.16)(40,-1.74)
\rput(125,18.75){$\Phi_n=\Phi_p$}
\psline[linestyle=dotted](130.92,37.11)(5,37.11)
\psline[linestyle=dotted](140.26,52.37)(5,52.37)
\psline{[-]}(5,37.11)(5,52.37)
\rput(21.25,45){$Q(\tilde\Gamma^\delta_{kjn})$}
\psline[linestyle=dotted](53.75,-5)(53.75,1.25)
\end{pspicture}
~\\
\vspace{10mm}

To come to a similar conclusion concerning $j$ and $j'$, we consider the projection on the $y-$axis, denoted $Q(x,y,z)=y$. Then - cf. \eqref{xsidef} -
$Q(\xi_{kjn})=\Phi_n(x_{kj})$ and $Q(\xi_{lip})=\Phi_p(x_{li})$, analogue for the primed coordinates.

\begin{stepnr}\label{step6}~
\begin{enumerate}
	\item $Q(\xi_{lip})-Q(\xi_{l'i'p}) \lesssim \delta 2^{l\mhalb}$
	\item $Q(\xi_{kj'n})-Q(\xi_{kjn}) \approx \delta 2^{k\mhalb} (j'-j)$   
\end{enumerate}\end{stepnr}
\underline{Part (ii):}
Again, we use the mean value theorem: there exists a $\tilde x$ between $x_{kj}$ and $x_{kj'}$ fulfilling
\begin{align*}
	Q(\xi_{kj'n})-Q(\xi_{kjn})
	= \Phi_n(x_{kj'})-\Phi_n(x_{kj}) 
	= \Phi_n'(\tilde x)(x_{kj'}-x_{kj})
	= \Phi'\left(\frac{\tilde x}{n\beta}\right)(x_{kj'}-x_{kj}).
\end{align*}
Notice that $\tilde x\approx x_k \approx 2^k\gamma$, i.e.
\begin{align*}
	Q(\xi_{kj'n})-Q(\xi_{kjn})
	\approx& (2^k\gamma)^{m-1}(x_{kj'}-x_{kj})	\\
	=& 2^{k(m-1)}\gamma^{m-1}(j'-j)2^{k(1-\mhalb)}\gamma	\\
	=& \delta 2^{k\mhalb}(j'-j).
\end{align*}
\underline{Part (i):}
In this case, the mean value theorem provides a $\bar x\approx x_l\approx 2^l\gamma$ such that
\begin{align*}
	Q(\xi_{lip})-Q(\xi_{l'i'p})
	=& \Phi_p(x_{li})-\Phi_p(x_{l'i'}) 	\\
	=& \Phi_p'(\bar x)(x_{li'}-x_{l'i'})	\\
	\approx& (2^l\gamma)^{m-1} ~ 2^{-l'(\mhalb-1)} \gamma \Delta_{ll'ii'}	\\
	=& \delta ~ 2^{l\mhalb} ~ 2^{(l-l')(\mhalb-1)} ~ \Delta_{ll'ii'}.		
\end{align*}
Since $\Delta_{ll'ii'}\lesssim1$ by Step \ref{step4} and $l'\approx l$ according to Step \ref{step5}, the equation implies
\begin{align*}
	Q(\xi_{lip})-Q(\xi_{l'i'p})
	\lesssim \delta 2^{l\mhalb},
\end{align*}
thus completing Step \ref{step6}.\\

\begin{stepnr}\label{step7}
$|j-j'|\lesssim 1$
\end{stepnr}
We have
\begin{eqnarray*}
	\diam_Q(\Gamma_{kjn}^\delta)	&=&	\sup_{\xi,\eta\in\Gamma_{kjn}^\delta}|Q(\xi)-Q(\eta)|	\\
	&=& Q(x_{kj},\Phi_n(x_{kj})+\landau(\delta),n\beta)-Q(x_{k,j-1},\Phi_n(x_{k,j-1}),n\beta) \\
	&=& \Phi_n(x_{kj}) + \landau(\delta) - \Phi_n(x_{k,j-1}) 	\\
	&\underset{Step~\ref{step6}(ii)}{\approx}& \delta 2^{k\mhalb} + \landau(\delta)\\
	&\lesssim& \delta 2^{k\mhalb}.
\end{eqnarray*}
Using Step \ref{step6}(i) we observe
$\diam_Q(\Gamma_{lip}^\delta)\lesssim \delta 2^{l\mhalb} \leq \delta 2^{k\mhalb}$ and
$\diam_Q(\Gamma_{l'i'p}^\delta)\lesssim \delta 2^{k\mhalb}$. Now we again apply Lemma \ref{überlapplemma}:
\begin{eqnarray*}
	Q(\xi_{kj'n})-Q(\xi_{kjn}) &+& Q(\xi_{l'i'p})-Q(\xi_{lip}) \\
	&\leq& diam_Q(\Gamma_{kj'n}^\delta)+diam_Q(\Gamma_{kjn}^\delta)+diam_Q(\Gamma_{l'i'p}^\delta)+diam_Q(\Gamma_{lip}^\delta)	\\
	&\lesssim& \delta 2^{k\mhalb}.
\end{eqnarray*}
It follows that
\begin{eqnarray*}
	\delta 2^{k\mhalb} (j'-j)
	&\underset{Step~\ref{step6}(ii)}{\approx}& Q(\xi_{kj'n})-Q(\xi_{kjn}) \\
	&\lesssim& Q(\xi_{lip})-Q(\xi_{l'i'p}) + \delta 2^{k\mhalb}	\\
	&\underset{Step~\ref{step6}(i)}{\lesssim}& \delta 2^{l\mhalb} + \delta 2^{k\mhalb}	\\
	&\lesssim& \delta 2^{k\mhalb},
\end{eqnarray*}
and thus $j'-j\lesssim 1$. This already implies the desired estimate, since $j\leq j'+1$ according to Step \ref{step1}.\\
%
Together with Step \ref{step5}, Step \ref{step7} completes the proof of Lemma \ref{ueberlapprest}.

\end{bew}


\subsection{Summary of the results}
\begin{kornr}
For all $\xi\in\R^3$, we have
\begin{align}
	\#\{(\alpha,\mu)|\xi\in\Gammaplus\} \leq 18 C \beta\inv,
\end{align}
where $C$ is the constant from Lemma \ref{ueberlapprest}.
\end{kornr}
\begin{bew}
For $\xi\in\R^3$ let
\begin{align*}
	M=&\{(n,p)|\exists k,l,j,i:\xi\in\Gammaplus\}, \\
	M(n,p)=&\{k|\exists l,j,i:\xi\in\Gammaplus,~l\leq k\}, \\  
	M(n,p;k)=&\{(l,j,i)|\xi\in\Gammaplus\}.
\end{align*}
Then the previous Lemma \ref{ueberlappz}, \ref{ueberlappk} and \ref{ueberlapprest}
states
\begin{align}
	\#M \leq& 3\beta\inv \\
	\#M(n,p) \leq& 3 \qquad\forall n,p\\
	\#M(n,p;k) \leq& C \qquad\forall n,p,k.
\end{align}
We proceed by
\begin{align*}
	\{(\alpha,\mu)|\xi\in\Gammaplus\} 
	=& \{(\alpha,\mu)|\xi\in\Gammaplus,~l\leq k\} \cup \{(\alpha,\mu)|\xi\in\Gammaplus,~l\geq k\} \\
	=& \{(\alpha,\mu)|\xi\in\Gammaplus,~l\leq k\} \cup \{(\alpha,\mu)|
												 \xi\in\Gamma^\delta_\mu+\Gamma^\delta_\alpha,~k\leq l\},										 
\end{align*}
and
\begin{align*}
	\{(\alpha,\mu)|\xi\in\Gammaplus,~l\leq k\}
	=& \{(\alpha,\mu)|(n,p)\in M,~k\in M(n,p),~(l,j,i)\in M(n,p;k) \}.
\end{align*}
It follows
\begin{align*}
	\#\{(\alpha,\mu)|\xi\in\Gammaplus\}
	\leq& 2\#\{(\alpha,\mu)|\xi\in\Gammaplus,~l\leq k\} \\
	=& 2 {\sum_{(n,p)\in M}}\, {\sum_{k\in M(n,p)}} \#M(n,p;k) \\
	\leq& 2 \cdot 3\beta\inv \cdot 3 \cdot C = 18C\beta\inv
\end{align*}
using Fubini's theorem.
\end{bew}

If we enlarge the sets $\Gamma_\alpha^\delta$ in the right manner, the statement essentially remains valid:

\begin{kornr}
For the sets $\Gamma_{kjn}^\delta$, we introduce their adjacent sets
$$\Gamma_{kjn}^\delta(u,v,w)=(0,v\delta,0)+\Gamma_{k,j+u,n+w}^\delta,\qquad u,v,w\in\{-1,0,1\},$$ as well as their "doubling" $$G_{kjn}^\delta=\bigcup\limits_{u,v,w\in\{-1,0,1\}} \Gamma_{kjn}^\delta(u,v,w).$$ 
Then for every $\xi\in\R^3$:
\begin{align*}
	\#\{(\alpha,\mu)|\xi\in G_\alpha^\delta+G_\mu^\delta\} \leq 3^3\cdot18 C \beta\inv.
\end{align*}
\end{kornr}

The last variation of this lemma is the version we want to apply finally.\\
\begin{kornr}\label{ueberlapp}
Let $1\leq s< \infty$ and $G_{kjn}^\delta$ like above. Then there exists a constant $C_s>0$ satisfying the following:
Let $f_{\alpha\mu}\in C_0^\infty(\R^3)$, $k,l=\krange$, $j=0,\ldots 2^{k\mhalb}-1$, $i=0,\ldots 2^{l\mhalb}-1$, $n,p=\zrange$, be non-negative functions fulfilling
\begin{align*}
	\supp f_{\alpha\mu}\subset G_\alpha^\delta+G_\mu^\delta 
\end{align*}
Then for all $\xi\in\R^3$
\begin{align*}
	\left( \sum_{\alpha\mu} f_{\alpha\mu}(\xi) \right)^s \leq C_s \beta^{1-s}\sum_{\alpha\mu} |f_{\alpha\mu}(\xi)|^s .
\end{align*}
\end{kornr}

\section{$L^p$-estimates for convolutions}\label{kapitel4}

Consider two cuboids in $\R^3$ with two short and one long edge. Both cuboids shall lie at parallel planes (we will concretise this soon). If we form the convolution of two functions, each one supported in one of the cuboids, what can we say about the $L^p$-Norm, depending on size and relative position?\\

For $i=1,2$ let $A_i=\{(x,y):|x|\leq\tfrac{\gamma_i}{2},~|y-m_ix|\leq\tfrac{\delta}{2}\}$. We assume that 
$\delta\ll\gamma_2,\gamma_1$ and moreover, that the slope $m_i$ of the boxes is bounded by an absolute constant (i.e. not depending on $\delta$ and $\gamma_i$). For convenience, let $\gamma_2\leq\gamma_1$.\\

\begin{center}\begin{picture}(8,8)(-4,-4)
	\refstepcounter{abb}\label{abbfaltung}
	\put(0.4,-2.9){\textbf{Figure \thesection.\theabb:} Parallelogram $A_i$}
	\newsavebox{\koord}
	\savebox{\koord}(0,3)[bl]{
	\put(-4,0){\vector(1,0){8}}
	\put(0,-3){\vector(0,1){7}}	 
	\put(-0.3,3.6){$y$} \put(3.7,-0.3){$x$}  }
	\put(0,-3){\usebox{\koord}}		\put(2,2.1){$A_i$}
	\put(-3,-2.5){\line(1,1){6}}\put(-3,-3.5){\line(1,1){6}}
	\put(3,2.5){\line(0,1){1}}\put(-3,-3.5){\line(0,1){1}}
	\put(-3,1){\line(1,0){6}}\put(-3,0.9){\line(0,1){0.2}}\put(3,0.9){\line(0,1){0.2}}
	\put(-2,1.2){$\gamma_i$}
	\put(-3.5,-3.5){\line(0,1){1}}\put(-3.6,-3.5){\line(1,0){0.2}}\put(-3.6,-2.5){\line(1,0){0.2}}
	\put(-3.8,-3.1){$\delta$}
\end{picture}\end{center}

Furthermore, we configure the position of the figures relative to each other. 
Let $\alpha=\angle(A_1,A_2)$ be the angle between the parallelograms, i.e. between their longer sides. 
The assumption of the boundedness of the slopes guarantees
\begin{align}\label{sinalfa}
	\sin\alpha\approx\alpha\approx\tan\alpha.
\end{align}

\stepcounter{abb}
\begin{center}\begin{picture}(8,8)(-4,-4)
\put(0.4,-2.9){\textbf{Figure \thesection.\theabb:} Parallelograms $A_i$}
\label{relativelage}
\newsavebox{\koorddd}
	\savebox{\koorddd}(0,3)[bl]{
	\put(-4,0){\vector(1,0){8}}
	\put(0,-3){\vector(0,1){7}}	 
	\put(-0.3,3.6){$y$} \put(3.7,-0.3){$x$}  }
\put(0,-3){\usebox{\koorddd}}
\put(-3,-2.75){\line(1,1){6}}\put(-3,-3.25){\line(1,1){6}}
\put(3,2.75){\line(0,1){0.5}}\put(-3,-3.25){\line(0,1){0.5}}
\put(1,2){$A_2$}
\put(-5,-1.25){\line(5,1){10}}\put(-5,-0.75){\line(5,1){10}}
\put(1,-0.6){$A_1$}
\qbezier(2,1.75)(2.375,1.375)(2.5,0.75)
\put(1.7,0.9){$\alpha$}
\end{picture}\end{center}

Now we introduce the parallelepipeds $Q_i=A_i\times(0,\beta)\subset\R^3$. Let $\xi_i\in\R^3$ and $\phi_i$ 
a bump function adopted to $\xi_i+Q_i$.

The main result of this section is the following:
\begin{lemnr}\label{faltungslemma}
	Let $1\leq s<\infty$. Then
	\begin{align}\label{faltung1}
		\int_{\R^3} |\phi_1\ast\phi_2|^s \dd x \lesssim 
		\frac{(\beta\delta)^{(s+1)}}{(1+ \tfrac{\gamma_2}{\delta}\alpha)^{s-1}} \gamma_2^s \gamma_1.
	\end{align}
\end{lemnr}

At the beginning a simple remark:		
\begin{bemnr}\label{falt1} 
	For $C>0$ and $i=1,2$ let $B_i$ be symmetric $(B_i=-B_i)$ subsets of $\R^n$, $x_i\in\R^n$, 		$\supp \psi_i\subset x_i+B_i$ and $\|\psi_i\|_\infty\leq C$. Then 
	\begin{align}\label{faltung100}
		\|\psi_1\ast\psi_2\|_\infty \leq C^2 \sup_z |B_1 \cap (z+B_2)|.
	\end{align}
\end{bemnr}
~\\


Further we need to estimate the size of the intersection $|A_1\cap A_2|$ and the sum $|A_1+A_2|$, which is by some basic geometric considerations.
\begin{lemnr}\label{falt2} We have
\begin{align*}
	|A_1\cap A_2| \lesssim \frac{\delta^2\gamma_2}{\max\{\delta,\gamma_2\alpha\}} 
		\approx \frac{\delta^2\gamma_2}{\delta+ \gamma_2\alpha}.
\end{align*}
\end{lemnr}

\begin{lemnr} \label{falt3} 
	$$|A_1+A_2| \lesssim  \gamma_1(\delta+\gamma_2\alpha).$$
\end{lemnr}

This Observations at hand prove Lemma \ref{faltungslemma}:
\begin{bew}
	\begin{eqnarray*}
					\int_{\R^3} |\phi_1\ast\phi_2|^s \dd x 
		&\leq&					\|\phi_1\ast\phi_2\|_\infty^s ~ |\supp(\phi_1\ast\phi_2)| \\
		&\underset{\text{Lemma }\ref{falt1}}{\lesssim}&	
							\Big(\sup_{\eta\in\R^2\x\R} |Q_1\cap(\eta+Q_2)|\Big)^s ~ |\xi_1+Q_1+\xi_2+Q_2| 	\\
		&=&	|Q_1\cap Q_2|^s ~ |Q_1+Q_2|		\\
		&=&	|\big(A_1\cap	A_2\big)\x(0,\beta)|^s ~ |\big(A_1\x(0,\beta)\big)+\big(A_2\x(0,\beta)\big)|\\
		&=&	|A_1\cap	A_2|^s\beta^s ~ |\big(A_1+A_2\big)\x(0,2\beta)|	\\
		&\underset{\text{Lemma }\ref{falt2}}{\lesssim}&	
						\bigg(\frac{\delta^2\gamma_2}{\delta+\gamma_2\alpha} \bigg)^s ~
						 |\big(A_1+A_2\big)|~\beta^{s+1}	\\
		&\underset{\text{Lemma }\ref{falt3}}{\lesssim}&	
						\bigg(\frac{\delta^2\gamma_2}{\delta+\gamma_2\alpha} \bigg)^s~
					 		\gamma_1(\delta+\gamma_2\alpha)~\beta^{s+1}\\
		&=&	\frac{\delta^{2s}\,\gamma_2^s\,\gamma_1}{(\delta+\gamma_2\alpha)^{s-1}}~\beta^{s+1} \\
		&=&  \frac{(\beta\delta)^{s+1}\,\gamma_2^s\,\gamma_1}{(1+\tfrac{\gamma_2}{\delta}\alpha)^{s-1}}.
		\end{eqnarray*}
\end{bew}

\section{Proof of the main theorem}
Before we start to complete the proof of Theorem \ref{superlemma}, we need a further lemma, which concretises the general results from the previous chapter in our special situation.

\subsection{Application of the results from Chapter \ref{kapitel4}}
\begin{lemnr}\label{korfuenf}
Let $1\leq s<\infty$. The functions $\phi_{kjn}$ introduced in the previous chapter satisfy the following estimates:
\begin{enumerate}\renewcommand{\labelenumi}{(\roman{enumi})}
\item If $|k-l|\gg1$, then
		\begin{align}\label{25oct1022i}
				\int|\phi_{kjn}\ast\phi_{lip}|^s\dd x  \lesssim&
				(\beta\delta\gamma)^{s+1} 2^{\frac{k+l}{2}(1+s-sm)} \  2^{-\frac{|k-l|}{2}(m-1)(s-1)}. 
		\end{align}
		
\item If $|k-l|\lesssim 1$, we find either a function $f:I_k\to I_l$ with $\forall i\in 		I_l:|\{j:f(j)=i\}|\lesssim 1$ or a function $g:I_l\to I_k$ with $\forall j\in I_k:|\{i:g(i)=j\}|\lesssim 1$ such that
\begin{align}\label{25oct1022ii}
	\int|\phi_{kjn}\ast\phi_{lip}|^s\dd x  \lesssim&
	(\beta\delta\gamma)^{s+1} 2^{\frac{k+l}{2}(1+s-sm)} 
								\ \left[\frac{2^{\frac{k+l}{2}\mhalb}}{1+h(i,j)}\right]^{s-1},
\end{align}
where 
\begin{align}
	h(i,j)=|f(j)-i| \qquad\text{or}\qquad  h(i,j)=|j-g(i)| .
\end{align}
\end{enumerate}
\end{lemnr}


\begin{center}\begin{picture}(8,8)(-4,-4)
\stepcounter{abb}
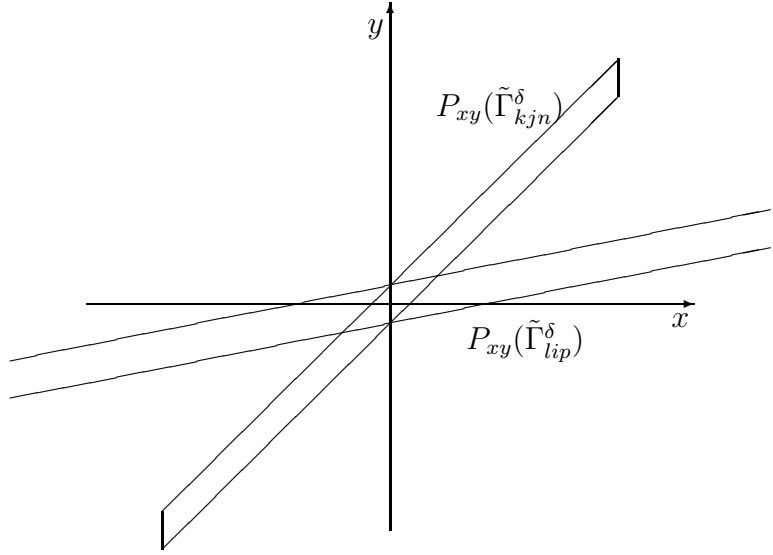
\put(-5.4,-3.75){\textbf{Figure \thesection.\theabb:} Angle between the translations to the origin of the projections}
\newsavebox{\koordd}
	\savebox{\koordd}(0,3)[bl]{
	\put(-4,0){\vector(1,0){8}}
	\put(0,-3){\vector(0,1){7}}	 
	\put(-0.3,3.6){$y$} \put(3.7,-0.3){$x$}  }
\put(0,-3){\usebox{\koordd}}
\put(-3,-2.75){\line(1,1){6}}\put(-3,-3.25){\line(1,1){6}}
\put(3,2.75){\line(0,1){0.5}}\put(-3,-3.25){\line(0,1){0.5}}
\put(0.6,2.5){$P_{xy}(\tilde\Gamma_{kjn}^\delta)$}
\put(-5,-1.25){\line(5,1){10}}\put(-5,-0.75){\line(5,1){10}}
\put(1,-0.6){$P_{xy}(\tilde\Gamma_{lip}^\delta)$}
\end{picture}\end{center}

\begin{bew} 
Recall the definition of $\tilde\Gamma_{kjn}^\delta=\{(x,y,z)\in\R^3|\, n\beta\leq z\leq(n+1)\beta,~	|y-\Phi_n(x)|\leq 10m\delta,~  x_{kj}\leq x\leq x_{k,j+1}\}$. 
The projections of these sets on $x$-$y-$space are contained in parallelograms with slope $\Phi'_n(x_{kj})$, thickness $\landau(\delta)$ and width $\gamma_{k}$. When we shift their centers to the origin, they intersect with angle
\begin{align}
	\alpha =& \angle(\tilde\Gamma_{kjn}^\delta,\tilde\Gamma_{lip}^\delta)
	= |\arctan \Phi'_n(x_{kj}) - \arctan \Phi'_p(x_{li})|.
\end{align}
Since $\Phi'$ and therefore as well $\Phi_n'=\Phi'\left(\frac{\cdot}{n\beta}\right)$ and $\Phi'_p$ are bounded from above and from below independently of $n$ and $p$, we have
\begin{align}\label{winkel}
	\alpha \approx& |\Phi'_n(x_{kj}) - \Phi'_p(x_{li})|.
\end{align}
At first we consider the case (i) $|k-l|\gg1$ and assume without loss of generality $l\leq k$, we obtain due to the dyadic nature of the construction
\begin{align}\label{11jan1247}
	\alpha \approx& \Phi'_n(x_{kj}) \approx (2^k\gamma)^{m-1} \nonumber\\
	=& \frac{\gamma^m}{2^{k(1-\mhalb)}\gamma} 2^{k\mhalb} \nonumber\\
	=& \frac{\delta}{\gamma_k} 2^{k\mhalb}.	
\end{align}

The case (ii) $|k-l|\lesssim 1$ is a bit more complicated. This is due to the fact that we decomposed
orthogonal to the $x-$axis, regardless the cone-like shape of the surface.
Formula \eqref{winkel} illustrates the difficulties: in the special case $n=p$ we get rid of the difference by the fundamental theorem of calculus, but unfortunately, the general case appears much harder.\\
Instead of this, we choose an other approach: let $k,l,n,p$ be fixed and $a_i=\Phi'_p(x_{li})$, $b_j=\Phi'_n(x_{kj})$. Then for all $j,j'$
\begin{eqnarray}\theoremsymbol{}
	|b_j-b_{j'}|&=&\left| \Phi'_n(x_{kj})	- \Phi'_n(x_{kj'}) \right| \nonumber\\
	&\approx& \Phi''(x_k) |x_{kj}-x_{kj'}| 			\nonumber\\
	&\approx& (2^k\gamma)^{m-2} \gamma_k |j-j'| \nonumber\\
	&=& 2^{k(\mhalb-1)}\gamma^{m-1} |j-j'| 			\nonumber\\
	&=& \frac{\delta}{\gamma_k} |j-j'|.
\end{eqnarray}
In the same manner we obtain
\begin{align}\label{11jan1311}
	|a_i-a_{i'}| \approx \frac{\delta}{\gamma_l} |i-i'| \overset{|k-l|\lesssim 1}{\approx} \frac{\delta}{\gamma_k} |i-i'|.
\end{align}
This basically means that we have a good idea how to compare the $a_i$ with each other, the same with the $b_j's$. But we lack control in comparing some $a_i$ with a $b_j$. Therefore
we apply the abstract result of Lemma \ref{folgenlemma}. At first we check the preconditions:

The sequences $a_i,i\in I_l$ and $b_j,j\in I_k$ are increasing, since (for instance)
$x_{li}<x_{li}+\gamma_{l}=x_{l,i+1}$ and $\Phi_p'$ is monotonously increasing. 
Thus we may apply Lemma \ref{folgenlemma} to the renormed sequences $\frac{\gamma_k}{\delta} a_i$ and $\frac{\gamma_k}{\delta} b_j$. Provided
$b_{2^{km/2}-1}-b_0\leq a_{2^{lm/2}-1}-a_0$ we get a function $f:I_k\to I_l$ almost injective, i.e. $\forall i\in I_l:|\{j\in I_k:\,f(j)=i\}|\lesssim 1$, and fulfilling
\begin{align}\label{11jan1313}
	|a_i-b_j|\geq \durch{2} |a_i-a_{f(j)}|\quad\forall i\in I_l\forall j\in I_k.
\end{align}
Consequently
\begin{eqnarray*}
	\alpha &\approx& |\Phi'_n(x_{kj}) - \Phi'_p(x_{li})| \\
	&=& |b_j-a_i| 	\nonumber\\
	&\overset{\eqref{11jan1313}}{\gtrsim}& |a_i-a_{f(j)}| 	\\
	&\overset{\eqref{11jan1311}}{\approx}& \frac{\delta}{\gamma_k} |i-f(j)|. 
\end{eqnarray*}
Provided $b_{2^{km/2}-1}-b_0> a_{2^{lm/2}-1}-a_0$ we switch the roles of $a_i$ and $b_j$ and obtain a correspondingly result; in any case, we obtain a function $h$ as desired and fulfilling
\begin{align}\label{21nov1336}
	\alpha \gtrsim \frac{\delta}{\gamma_k} h(i,j).
\end{align}
Combining \eqref{11jan1247} and \eqref{21nov1336} we see that
\begin{align}\label{11jan1510}
	\frac{\gamma_k}{\delta} \alpha 
	= \begin{cases} 2^{k\mhalb}, & \text{ if }l\ll k \\ 
						h(i,j), & \text{ if } |k-l|\lesssim 1 .
\end{cases}
\end{align}
Using Lemma \ref{faltungslemma} and taking into account $l\leq k$, i.e. $\gamma_k\leq\gamma_l$, we conclude 
\begin{align*}
	\int|\phi_{kjn}\ast\phi_{lip}|^s\dd x \lesssim&
		\frac{(\beta\delta)^{s+1}}{(1+\tfrac{\gamma_k}{\delta}\alpha)^{s-1}} \gamma_{k}^s \gamma_{l}  \\
	=& \frac{(\beta\delta\gamma)^{s+1}}{(1+\tfrac{\gamma_k}{\delta}\alpha)^{s-1}}
					 2^{k(1-\mhalb)s}2^{l(1-\mhalb)}  \\
	=& \frac{(\beta\delta\gamma)^{s+1}}{(1+\tfrac{\gamma_k}{\delta}\alpha)^{s-1}} 
										2^{\tfrac{k}{2}(2s-ms)}2^{\tfrac{l}{2}(2-m)}  \\
	=& (\beta\delta\gamma)^{s+1} 2^{\frac{k+l}{2}(1+s-sm)}
	\ \cdot\ \frac{2^{\tfrac{k}{2}(s-1)+\tfrac{l}{2}(ms-m-s+1)}}{(1+\tfrac{\gamma_k}{\delta}\alpha)^{s-1}}. 
\end{align*}
It remains to consider the second expression. In the case $l\ll k$, it is transformed by \eqref{11jan1510} into
\begin{align*}
	\frac{2^{\tfrac{k}{2}(s-1)+\tfrac{l}{2}(ms-m-s+1)}}{(1+\tfrac{\gamma_k}{\delta}\alpha)^{s-1}}
	\approx& \frac{2^{\tfrac{k}{2}(s-1)+\tfrac{l}{2}(ms-m-s+1)}}{2^{k\mhalb(s-1)}} \\
	=& 2^{\tfrac{k}{2}(s-1-ms+m)+\tfrac{l}{2}(ms-m-s+1)} \\
	=& 2^{-\tfrac{k-l}{2}(m-1)(s-1)},
\end{align*}
whereas in the case $|k-l|\lesssim1$, it gives
\begin{align*}
	\frac{2^{\tfrac{k}{2}(s-1)+\tfrac{l}{2}(ms-m-s+1)}}{(1+\tfrac{\gamma_k}{\delta}\alpha)^{s-1}}
	=& \frac{2^{\tfrac{k}{2}(s-1)+\tfrac{l}{2}(ms-m-s+1)}}{(1+h(i,j))^{s-1}} \\
	\overset{k\approx l}{\approx}& \frac{2^{\tfrac{k+l}{4}(s-1)+\tfrac{k+l}{4}(ms-m-s+1)}}{(1+h(i,j))^{s-1}} \\
	=& 	\frac{2^{\tfrac{k+l}{4}(ms-m)}}{(1+h(i,j))^{s-1}} \\
	=& 	\left( \frac{2^{\tfrac{k+l}{2}\mhalb}}{1+h(i,j)} \right)^{s-1},
\end{align*}
whereby the claim is verified.
\end{bew}

\subsection{Completing the proof}
One further intermediate step will be helpful since we get rid of the split-up in the two different cases from the previous lemma.
\begin{lemnr}\label{vereinigung}
Let $\frac{3}{p'}<\tfrac{1}{q}< \tfrac{1}{2}+\tfrac{1}{p'}$ and 
\begin{align}\label{sdefi}
	\frac{1}{s}+\frac{2}{p'}=1. 
\end{align}
Then for all $k,l,n,p$ and for all finite sequences $a\in\R^{I_k},b\in\R^{I_l}$ we have
\begin{align*}
	&\sum\limits_{ij} |a_j|^s|b_i|^s \int|\phi_{kjn}\ast\phi_{lip}|^s\dd x \\
	\lesssim& (\beta\delta\gamma)^{s+1} \|a\|_{q'}^s\|b\|_{q'}^s
	 2^{\frac{k+l}{2}(1+s-sm\,+\,m\left(1-\frac{s}{q'}\right))}
	 2^{-\frac{|k-l|}{2}(m-1)(s-1)}		.	
\end{align*}
\end{lemnr}
\begin{bew}
At first we consider the case $|k-l|\gg1$. Hölder's inequality implies
\begin{align*}
	\sum\limits_{j} |a_{j}|^s \leq&
	\left(\sum_{j} |a_{j}|^{q'}\right)^{\frac{s}{q'}}		\left(\sum_{j=0}^{2^{k\frac{m}{2}}-1}
			 \right)^{1-\frac{s}{q'}} \\
	=& \|a\|_{q'}^s 2^{\frac{k}{2}m\left(1-\frac{s}{q'}\right)},
\end{align*}
which, together with Lemma \ref{korfuenf}, results in the desired estimate.\\

In the case $|k-l|\leq 1$ we apply Lemma \ref{fracint} with parameter 
\begin{align}\label{rdefi}
	r=\frac{q'}{s}.
\end{align}
Since we assumed $\durch{q}<\durch{2}+\durch{p'}$, it follows
\begin{align}
	\durch{q'}=1-\durch{q}>\durch{2}-\durch{p'}=\durch{2}\left(1-\frac{2}{p'}\right)
	\overset{\eqref{sdefi}}{=}\durch{2s},
\end{align}
thus
\begin{align}\label{neunundfufzig}
	r=\frac{q'}{s}<2.
\end{align}
Furthermore, we have
\begin{eqnarray*}
	\durch{r'}&=& 1-\frac{s}{q'} = 1-s+\frac{s}{q} \\
	&\overset{ }{>}& 1-s+s\frac{3}{p'} = 1-s+\frac{3s}{2}\frac{2}{p'} \\
	&\overset{\eqref{sdefi}}{=}& 1-s+\frac{3s}{2}\left(1-\durch{s}\right) =
								 -(s-1)+\frac{3}{2}(s-1) \\
	&=&\frac{s-1}{2},
\end{eqnarray*}
i.e.
\begin{align}\label{Qdefi}
	Q:=(s-1)\frac{r'}{2}<1.
\end{align}
Assume for simplicity once more that $l\leq k$.
To apply the lemma, we have to analyse the kernel $G_{\frac{k+l}{2}}(x)=\frac{2^{\frac{k+l}{2}\mhalb}}{1+x}$ from \eqref{25oct1022ii}, $x\in[0,2^{k\mhalb}-1]$. It fulfills
\begin{align*}
	\int_0^{2^{k\mhalb}-1} G_{\frac{k+l}{2}}^Q\dd x \leq& \int_1^{2^{k\mhalb}} \left( \frac{2^{k\mhalb}}{x}\right)^Q \dd x 
	= 2^{k\mhalb} \int_{2^{-k\mhalb}}^1 x^{-Q} \dd x	\\
	\leq& 2^{k\mhalb}  \left[ \frac{x^{1-Q}}{1-Q} \right]_0^{1}  
	\approx 2^{k\mhalb} \overset{|k-l|\lesssim 1}{\approx} 2^{\frac{k+l}{2}\mhalb}, 
\end{align*}
	i.e.
\begin{align}\label{23nov1302}
	\| G_\frac{k+l}{2}^{s-1} \|_{\frac{r'}{2}} = \|G_\frac{k+l}{2}\|_Q^{s-1} 
		\lesssim 2^{\frac{k+l}{2}\mhalb\frac{s-1}{Q}} 
		\overset{\eqref{Qdefi}}{\leq} 2^{\frac{k+l}{2}\frac{m}{r'}}\overset{\eqref{rdefi}}{=}2^{\frac{k+l}{2} m\left(1-\frac{s}{q'}\right)}.
\end{align}
Next we again apply Lemma \ref{korfuenf}, but this time part (ii), though we just discuss the first case $h(i,j)=|f(j)-i|$. The lemma states
\begin{align}\label{27nov1632}
	&\sum\limits_{ij} |a_j|^s|b_i|^s \int|\phi_{kjn}\ast\phi_{lip}|^s\dd x \nonumber\\
	\lesssim& 	(\beta\delta\gamma)^{s+1} \sum\limits_{ij} |a_j|^s|b_i|^s 2^{\frac{k+l}{2}(1+s-sm)} 
								\left[\frac{2^{\frac{k+l}{2}\mhalb}}{1+|f(j)-i|}\right]^{s-1} 
\end{align}
and according to Lemma \ref{fracint} 
\begin{eqnarray*}
	\sum\limits_{ij} |a_j|^s|b_i|^s \left[\frac{2^{\frac{k+l}{2}\mhalb}}{1+|f(j)-i|}\right]^{s-1}
	&\lesssim& \|a^s\|_{r}\|b^s\|_r \|G_{\frac{k+l}{2}}^{s-1}\|_{\frac{r'}{2}}  \\
	&\overset{\eqref{23nov1302}}{\lesssim}&\|a^s\|_{r}\|b^s\|_r 2^{\frac{k+l}{2} m\left(1-\frac{s}{q'}\right)}\\
	&\overset{\eqref{rdefi}}{=}&\|a\|_{q'}^s\|b\|_{q'}^s 2^{\frac{k+l}{2} m\left(1-\frac{s}{q'}\right)}.
\end{eqnarray*}
Putting this into \eqref{27nov1632} yields the claim, when we insert the expression
$2^{-\frac{|k-l|}{2}(m-1)(s-1)}$, which is of constant order in case of $|k-l|\lesssim1$.
\end{bew}


~\\

Eventually we can complete the proof of Theorem \ref{superlemma}. Let us recall the statement:\\
\begin{satz}
Let $p'>m+1$, $\frac{3}{p'}<\tfrac{1}{q}< \tfrac{1}{2}+\tfrac{1}{p'}$, then
\begin{align}\label{31oct1215}
		\|\sum_{\alpha} a_\alpha \hat\phi_\alpha\|_{p'}	\lesssim&  
		\delta^{\durch{q}}\|\sum_{\alpha=(k,j,n)} 
							a_\alpha (\zweik\gamma)^{\durch{q}-\frac{m+1}{p'}} \phi_\alpha\|_{q'}
\end{align}
holds for all $\delta>0$ and for all sequences $a_\alpha$.
\end{satz}
\begin{bew}
Since $2<p'<\infty$, we find $s\in(1,\infty)$ such that 
\begin{align}\label{sdefi2}
	\frac{1}{s}+\frac{2}{p'}=1.
\end{align}

Since even $4<p'$, $\frac{p'}{2}>2$ holds, thus we can apply Youngs inequality:
\begin{eqnarray}
	\|\sum_{\alpha} a_\alpha \hat\phi_\alpha\|_{p'}^{2s}
	&=& \|\sum_{\alpha,\mu} a_\alpha a_\mu \hat\phi_\alpha \hat\phi_\mu \|_{\tfrac{p'}{2}}^{s} \nonumber\\
	&=& \|\FT\Big( \sum_{\alpha,\mu} a_\alpha a_\mu \phi_\alpha\ast\phi_\mu\Big) \|_{\tfrac{p'}{2}}^{s} \nonumber\\
	&\overset{Young}{\leq}& \|\sum_{\alpha,\mu} a_\alpha a_\mu \phi_\alpha\ast\phi_\mu \|_{s}^{s} \nonumber\\
	&=& \int\Big| \sum_{\alpha,\mu} a_\alpha a_\mu \phi_\alpha\ast\phi_\mu \Big|^s\dd x. \nonumber
\end{eqnarray}	
Now we exploit the estimate from Chapter \ref{kap:ueberlapp} of the overlap of the supports of the functions $\phi_\alpha\ast\phi_\mu$:
\begin{eqnarray}
\|\sum_{\alpha} a_\alpha \hat\phi_\alpha\|_{p'}^{2s}
	&\leq& \int\Big| \sum_{\alpha,\mu} a_\alpha a_\mu \phi_\alpha\ast\phi_\mu \Big|^s\dd x. \nonumber \\
	&\overset{\text{Corollary }\ref{ueberlapp}}	{\lesssim}& \beta^{1-s}
	\sum_{\alpha,\mu} |a_\alpha|^s |a_\mu|^s  \int\left| \phi_\alpha\ast\phi_\mu \right|^s\dd x \nonumber\\
	&\overset{\text{Lemma }\ref{vereinigung}}{\lesssim}& \beta^{2}(\delta\gamma)^{s+1}
	\sum_{klnp} \|a_{kn}\|_{q'}^s \|a_{lp}\|_{q'}^s 2^{\frac{k+l}{2}(1+s-sm\,+ 
					\,m\left(1-\frac{s}{q'}\right))} 2^{-\frac{|k-l|}{2}(m-1)(s-1)}, \nonumber
\end{eqnarray}
where $a_{kn}=(a_{kjn})_{j\in I_k}$, $a_{lp}=(a_{lip})_{i\in I_l}$. 
Notice that for fixed $k$ Hölder's inequality implies
\begin{align}\label{zhoelder}
	\sum\limits_n \|a_{kn}\|_{q'}^s \leq&
	\left(\sum_{n} \|a_{kn}\|_{q'}^{q'}\right)^{\frac{s}{q'}}	
			\left(\sum_{n=\durch{\beta}}^{\frac{2}{\beta}-1} 1 \right)^{1-\frac{s}{q'}} \nonumber\\
	=& \|a_{k}\|_{q'}^s \left(\durch{\beta}\right)^{1-\frac{s}{q'}} 
	=\|a_{k}\|_{q'}^s \beta^{\frac{s}{q'}-1},
\end{align}
where $a_{k}=(a_{kjn})_{j,n}$. 
Implementing this in our equations yields
\begin{align}
	\|\sum_{\alpha} a_\alpha \hat\phi_\alpha\|_{p'}^{2s} \lesssim&
	\beta^{2(\frac{s}{q'}-1)} \beta^2(\delta\gamma)^{(s+1)}\sum\limits_{k,l}\|a_k\|_{q'}^s \|a_l\|_{q'}^s
	2^{\frac{k+l}{2}(1+s-sm+m\left(1-\frac{s}{q'}\right))} 2^{-|k-l|\frac{(m-1)(s-1)}{2}} \nonumber\\
	=&\beta^{\frac{2s}{q'}} (\delta\gamma)^{(s+1)} \sum\limits_{k,l} \|a_k\|_{q'}^s \|a_l\|_{q'}^s
	2^{\frac{k+l}{2}[s(1-m)+m+1-m\frac{s}{q'}]} 2^{-|k-l|\frac{(m-1)(s-1)}{2}} \nonumber\\
	=&\beta^{\frac{2s}{q'}} (\delta\gamma)^{(s+1)} \sum\limits_{k,l} 
	\|a_k 2^{k[\frac{1-m}{2}+\frac{m+1}{2s}-\frac{m}{2q'}]}\|_{q'}^s 
	\|a_l 2^{l[\frac{1-m}{2}+\frac{m+1}{2s}-\frac{m}{2q'}]}\|_{q'}^s 
	 2^{-|k-l|\frac{(m-1)(s-1)}{2}} .
\end{align}
We again apply Lemma \ref{fracint} for $r=\frac{q'}{s}$ (recall $r<2$, cf. \eqref{neunundfufzig}):
\begin{align}\label{31oct0957}
	\|\sum_{\alpha} a_\alpha \hat\phi_\alpha\|_{p'}^{2s} \lesssim\quad&
	\beta^{\frac{2s}{q'}} (\delta\gamma)^{s+1} \left( \sum\limits_k \|a_k 2^{k[\frac{1-m}{2} +\frac{m+1}{2s}-\frac{m}{2q'}]} \|_{q'}^{q'} \right)^{\frac{2s}{q'}} \nonumber\\
	{=}\quad&
	 (\delta\gamma)^{s+1-\frac{2s}{q'}}
	 \gamma^{-\left(\frac{1-m}{2}+\frac{m+1}{2s}-\frac{1}{q'}\right)2s}
	 \left( \sum\limits_k
	  \|a_k (2^k\gamma)^{\frac{1-m}{2}+\frac{m+1}{2s}-\frac{1}{q'}} \|_{q'}^{q'} 
	   2^{k(1-\frac{m}{2})}\gamma\beta\delta\right)^{\frac{2s}{q'}} \nonumber\\
	 \overset{\ \ \delta=\gamma^m}{=}&
	 \gamma^{(m+1)\left(s+1-\frac{2s}{q'}\right)-
	 		\left(\frac{1-m}{2}+\frac{m+1}{2s}-\frac{1}{q'}\right)2s} 
	 \left( \sum\limits_k \|a_k (2^k\gamma)^{\frac{1-m}{2}+\frac{m+1}{2s}-\frac{1}{q'}} \|_{q'}^{q'}  2^{k(1-\frac{m}{2})}\gamma\beta\delta\right)^{\frac{2s}{q'}}. 
\end{align}
We determine the exponents: it holds
\begin{align}
	&(m+1)\left(s+1-\frac{2s}{q'}\right)-\left(\frac{1-m}{2}+\frac{m+1}{2s}-\frac{1}{q'}\right)2s
									 \nonumber\\
	=&2s\left[ (m+1)\left(\durch{2}+\frac{1}{2s}-\frac{1}{q'}\right)
											-\frac{1-m}{2}-\frac{m+1}{2s}+\frac{1}{q'} \right] \nonumber\\
	=& 2s\left[ \frac{m+1}{2}-\frac{m}{q'}-\frac{1-m}{2} \right] \nonumber\\
	=& 2s\,m\left[1-\durch{q'}\right]=\frac{2sm}{q},	
\end{align}
and
\begin{align}
	\frac{1-m}{2}+\frac{m+1}{2s}-\frac{1}{q'}=\frac{1}{q}-\frac{m+1}{2}+\frac{m+1}{2s}
	=\frac{1}{q}-\frac{m+1}{2}\left(1-\durch{s}\right) \overset{\eqref{sdefi2}}{=} 
		\frac{1}{q}-\frac{m+1}{p'}.
\end{align}
A further computation shows
\begin{align}
	\int \phi_\alpha^{q'}\dd x \approx |\Gamma_\alpha^\delta|=|\Gamma_{kjn}^\delta| \approx \gamma_k\beta\delta = 2^{k(1-\frac{m}{2})}\gamma\beta\delta,
\end{align}
thus \eqref{31oct0957} translates into
\begin{align*}
	\|\sum_{\alpha} a_\alpha \hat\phi_\alpha\|_{p'}^{2s} \lesssim&
	\gamma^{m\frac{2s}{q}} 
	 \left( \sum\limits_k \|a_k (2^k\gamma)^{\frac{1}{q}-\frac{m+1}{p'}} \|_{q'}^{q'}
	 			  \int \phi_\alpha^{q'}\dd x  \right)^{\frac{2s}{q'}} \\
	=&\delta^{\frac{2s}{q}} 
	 \left( \int \sum\limits_{\alpha}  \left|a_\alpha (2^k\gamma)^{\frac{1}{q}-\frac{m+1}{p'}} 
	 			  \phi_\alpha\right|^{q'}    \dd x  \right)^{\frac{2s}{q'}} \\
	\approx& \left(\delta^{\durch{q}} \|\sum_\alpha a_\alpha (\zweik\gamma)^{\frac{1}{q}-\frac{m+1}{p'}} \phi_\alpha \|_{q'}\right)^{2s},
\end{align*}
completing the proof of Theorem \ref{superlemma}.
\end{bew}

\section{On the necessity of the Lorentz space}
\label{kegelnotwendig}

We will give a short explanation why the weak-type estimate in Theorem \ref{thm0} for certain values of $p$ and $q$ cannot be improved to the strong type estimate.

\begin{lemnr}\label{7murch1654}
Let $m\geq1$, $s>0$, $0\leq\alpha<1$ and let $R\gg1$ be sufficiently large. Then
\begin{enumerate}
	\item \begin{align}
		\left| \Re \int_{\frac{1}{R}}^R e^{-ix^m}x^{-\alpha} \dd x\right| \geq c >0.
				\end{align}	
			and	\begin{align}
		\left| \int_{\frac{1}{R}}^R e^{-ix^m}x^{-\alpha} \dd x\right| \leq C.
				\end{align}	
	
	\item For $1\ll u^m\ll v$ we have 
			\begin{align}
			\left| \Re  \int_{\frac{1}{R}}^R e^{-i(x^m-uv^{-\frac{1}{m}}x)}x^{-\alpha} \left(1-\frac{m\log x}{\log v}\right)^{-s} \dd x\right| \geq \frac{c}{2}.
				\end{align}	
			and \begin{align}
			\left| \int_{\frac{1}{R}}^R e^{-i(x^m-uv^{-\frac{1}{m}}x)}x^{-\alpha} \left(1-\frac{m\log x}{\log v}\right)^{-s} \dd x\right| \leq 2C.
				\end{align}	
\end{enumerate}
\end{lemnr}

\begin{bew}
	To prove (i), introduce $\beta=\frac{\alpha-1}{m}+1\in(0,1)$. Using a contour integral, we can transform
\begin{align*}
	m \int_{R^{-\frac{1}{m}}}^{R^\frac{1}{m}} e^{-ix^m}x^{-\alpha} \dd x
	= \int_{\frac{1}{R}}^R e^{-iy}y^{-\beta} \dd y
\end{align*}
into
\begin{align}\label{7murch1556}
	\int_{\frac{1}{R}}^R e^{-t} t^{-\beta}\dd t.
\end{align}
Observe that
\begin{align}
	\int_{0}^\infty e^{-t} t^{-\beta}\dd t = \gamma(1-\beta)\neq 0,
\end{align}
whereas the remaining contributions are small as $R\to\infty$.\\
It is easy to deduce (ii) from (i): If $R$ is fixed, we may choose $\frac{u}{v^\frac{1}{m}}$ small enough and $v$ large enough such that the functions under the integrals in (i) and (ii) differ at most
$\landau(R^{-1})$.~
\end{bew}

\begin{lemnr}\label{6murch1642}
	Let $m\geq2$, $s\geq 0$, $0\leq\alpha<1$ and $1\ll u^m\ll v$. Then
	\begin{align}
		\left| \Re \int_0^{\frac{1}{2}} e^{i(ux-vx^m)}x^{-\alpha} |\log(x)|^{-s} \dd x\right|
		\gtrsim v^{\frac{\alpha-1}{m}} \log^{-s}(v)
	\end{align}
and
	\begin{align}
		\left| \int_0^{\frac{1}{2}} e^{i(ux-vx^m)}x^{-\alpha} |\log(x)|^{-s} \dd x\right|
		\lesssim v^{\frac{\alpha-1}{m}} \log^{-s}(v).
	\end{align}
\end{lemnr}
\begin{bew}
Applying the change of variables $x=v^{-\frac{1}{m}}y$ we see that
\begin{align*}
	\int_0^{\frac{1}{2}} e^{i(ux-vx^m)}x^{-\alpha} |\log(x)|^{-s} \dd x
	\approx& v^{\frac{\alpha-1}{m}} \log^{-s}(v)
			\int_0^{\frac{1}{2}v^{\frac{1}{m}}} e^{i\left(uv^{-\frac{1}{m}}y-y^m\right)} y^{-\alpha} 
																				\left|1-\frac{m\log(y)}{\log(v)}\right|^{-s} \dd y.
\end{align*}
We decompose the set of integration into $[0,1/R]$, $[1/R,R]$, $[R,v^\frac{1}{2m}]$ and $[v^\frac{1}{2m},\frac{1}{2}v^{\frac{1}{m}}]$.
We saw in Lemma \ref{7murch1654} that 
\begin{align}
	\left| \Re \int_\frac{1}{R}^R e^{i\left(uv^{-\frac{1}{m}}y-y^m\right)} y^{-\alpha} 
																				\left|1-\frac{m\log(y)}{\log(v)}\right|^{-s} \dd y \right|
	\gtrsim 1
\end{align}
and
\begin{align}
	\left| \int_\frac{1}{R}^R e^{i\left(uv^{-\frac{1}{m}}y-y^m\right)} y^{-\alpha} 
																				\left|1-\frac{m\log(y)}{\log(v)}\right|^{-s} \dd y \right|
	\lesssim 1. 
\end{align}

It remains to show that the other contributions to the integral are sufficiently small. 

For the first part of the integral, we have
\begin{align}
	\left| \int_0^\frac{1}{R} e^{i\left(uv^{-\frac{1}{m}}y-y^m\right)} y^{-\alpha} 
																				\left|1-\frac{m\log(y)}{\log(v)}\right|^{-s} \dd y \right|
	\lesssim \int_0^\frac{1}{R} y^{-\alpha} \dd y 
	\approx R^{\alpha-1} 
	\stackrel{R\to\infty}{\longrightarrow} 0.																				
\end{align}
For the third and fourth part where $y\geq R$ we have $$\frac{\dd}{\dd y} \left(y^m-uv^{-\frac{1}{m}}y\right) = my^{m-1}-uv^{-\frac{1}{m}} \gg R$$ and we may apply integration by parts. For any $c,d$ with $R\leq c\leq y\leq d\leq \frac{1}{2} v^\frac{1}{m}$ we have that
$1-\frac{m\log(y)}{\log(v)}\geq 1-\frac{m\log(d)}{\log(v)}$ and $\log v- m\log y\geq \log v- m\log\left(\frac{1}{2} v^\frac{1}{m}\right) = m\log 2$, thus
\begin{align*}
	\left| \int_c^d e^{i\left(uv^{-\frac{1}{m}}y-y^m\right)} y^{-\alpha} 
			\left|1-\frac{m\log(y)}{\log(v)}\right|^{-s} \dd y \right|
	\lesssim& \frac{1}{R} \left|1-\frac{m\log(d)}{\log(v)}\right|^{-s}	c^{-\alpha}.
\end{align*}
If we choose $c=R$ and $d=v^\frac{1}{2m}$, then 
$1-\frac{m\log(d)}{\log(v)} = \frac{1}{2}$. Hence we have
\begin{align}
	\left| \int_R^{v^\frac{1}{2m}} e^{i\left(uv^{-\frac{1}{m}}y-y^m\right)} y^{-\alpha} 
																				\left|1-\frac{m\log(y)}{\log(v)}\right|^{-s} \dd y \right|
	\lesssim R^{-1-\alpha}.
\end{align}
If we choose $c=v^\frac{1}{2m}$ and $d=\frac{1}{2}v^\frac{1}{m}$, then
$1-\frac{m\log(d)}{\log(v)} = \frac{m\log2}{\log v}$. We obtain
\begin{align}
	\left| \int_{v^\frac{1}{2m}}^{\frac{1}{2}v^\frac{1}{m}} 
	e^{i\left(uv^{-\frac{1}{m}}y-y^m\right)} y^{-\alpha} 
																					\left|1-\frac{m\log(y)}{\log(v)}\right|^{-s} \dd y \right|
	\lesssim \frac{1}{R} v^{-\frac{\alpha}{2m}} \log^{s}( v) \ll \frac{1}{R}.
\end{align}
~\end{bew}


\begin{lemnr}
	Let $\sigma$ be the Lebesgue measure of $S=\{(x,x^m z^{1-m},z)|x\in[0,1],\ z\in[1,2]\}$, and let $\frac{1}{q}=\frac{m+1}{p'}$. If the adjoint restriction estimate 
\begin{align}
	\|\widehat{f\dd\sigma}\|_{p'} \leq C \|f\|_{L^{q'}(S)}
\end{align}
holds true for some constant $C>0$ and any $f\in L^{q'}(S)$, then $q\geq p$.
\end{lemnr}

\begin{bew}
We will show that if $q<p$, then the restriction estimate fails. If $q<p$, i.e. $\frac{r}{q'}<\frac{1}{p'}\leq1$, we find $r>1$ such that
\begin{align}
	 s:=\frac{r}{q'}<\frac{1}{p'}.
\end{align}
Let $f(x,z)=f_1(x/z)\chi_{\{x/z\leq 1/2\}}$, $f_1(x)=x^{-\frac{1}{q'}} \log^{-s}(1/x)$, $0<x<\frac{1}{2}$. We claim that $f\in L^{q'}(S)$.\\
In the case $\frac{1}{q'}>0$, we have $s=\frac{r}{q'}>\frac{1}{q'}$, i.e. $1-sq'<0$. Thus we conclude
\begin{align*}
	 & \int_1^2 \int_0^1 |f(x,z)|^{q'} \dd x \dd z	\\
	\approx& \int_0^{\frac{1}{2}} x^{-1}\log^{-sq'}(1/x) \dd x \int_1^2\dd z	\\
	=& \int_2^{\infty} x^{-1}\log^{-sq'}(x) \dd x	\\
	\approx& \left[\log^{1-sq'}(x)\right]_\infty^2 < \infty.
\end{align*}
In the case $\frac{1}{q'}=0$, we have $s=0$ and $f_1\equiv 1$, i.e. 
$f=\chi_{\{x/z\leq 1/2\}} \in L^\infty(S)$.\\
Further 
we have
\begin{align*}
	 &	\widehat{f\dd\sigma}(u,-v,w)	\\
	=& \int_1^2 \int_0^{z/2} e^{i(ux+wz-vx^mz^{1-m})} f_1(x/z) \dd x\dd z	\\
	=& \int_1^2 z e^{izw} \int_0^{1/2} e^{iz(ux-vx^m)} f_1(x) \dd x\dd z	\\
	=& \int_1^2 z e^{izw} F_{u,v}(z) \dd z
\end{align*}
if we define 
\begin{align}
	F_{u,v}(z)= \int_0^{1/2} e^{iz(ux-vx^m)} f_1(x) \dd x = F_{zu,zv}(1).
\end{align}
Since $z\to \Re z e^{izw} F_{u,v}(z)$ is a continuous real-valued function, there exists $z_0\in[1,2]$ such that
\begin{align*}
	\Re \widehat{f\dd\sigma}(u,-v,w)
	= \Re \int_1^2 z e^{izw} F_{u,v}(z) \dd z
	= \Re z_0 e^{iz_0w} F_{u,v}(z_0).
\end{align*}
We know from Lemma \ref{6murch1642} that there exist constants $C_0,c_0>0$ such that for all $1\ll u^m \ll v$ 
\begin{align}
	|\Re F_{u,v}(1)| \geq& c_0 v^{\frac{1/q' -1}{m}} \log^{-s}(v)	\\
	| F_{u,v}(1)| \leq& C_0 v^{\frac{1/q' -1}{m}} \log^{-s}(v)	.
\end{align}
Observe that $z_0\in [1,2]$.
We conclude that for all $1\ll u^m \ll v$, $w\in[0,\e]$ we have
\begin{align*}
	|\Re \widehat{f\dd\sigma}(u,-v,w)|
	=& |\Re z_0 e^{iz_0w} F_{z_0u,z_0v}(1)|			\\
	\geq& z_0 \big[|\Re F_{z_0u,z_0v}(1)|\cos(z_0w)-|F_{z_0u,z_0v}(1)|\,|\sin(z_0w)| \big]	\\ 
	\geq& v^{\frac{1/q' -1}{m}} \log^{-s}(v) [c_0\cos(2\e)-2C_0\e]
	\gtrsim v^{\frac{1/q' -1}{m}} \log^{-s}(v)
\end{align*}
provided that $\e$ is small enough. Thus
\begin{align*}
	\|\widehat{f\dd\sigma}\|_{p'}^{p'}
	\geq&\int_0^\e\int_C^\infty\int_{1\ll u^m\ll v} |\widehat{f\dd\sigma}(u,-v,w)|^{p'} \dd u\dd v\dd w	\\
	\gtrsim& \e \int_C^\infty v^\frac{1}{m} v^{\frac{1/q' -1}{m}p'} \log^{-sp'}(v) \dd v.
\end{align*}
Observe 
\begin{align}
	\frac{1}{m} + \frac{1/q' -1}{m}p' 
	= \frac{1}{m} - \frac{p'}{qm} 
	= \frac{1}{m} - \frac{m+1}{m} 
	= -1.
\end{align}
We conclude
\begin{align}
	\|\widehat{f\dd\sigma}\|_{p'}^{p'} \gtrsim \e \int_C^\infty v^{-1} \log^{-sp'}(v)\dd v 
	\approx \e \left[\log^{1-sp'}(v)\right]_C^\infty = \infty
\end{align}
since $1-sp'>0$.~
\end{bew}

\section{Appendix}

The first lemma is a slightly modification of a well known fact. Nevertheless we will go through the simple proof. 

\begin{lemnr}\label{fracint}
Let $1\leq r<2$, and $f:\N\to\N$ "almost" injective, i.e. there exists a constant $C>0$ such that
\begin{align}\label{fastinjektiv}
	\forall l:|\{k\in\N:f(k)=l\}|\leq C.
\end{align}
Then
\begin{align*}
	\left| \sum_{k,l} a_k G_{f(k)-l} b_l \right| \leq C^{\frac{1}{r'}} \|a\|_r \|G\|_{\tfrac{r'}{2},\infty} \|b\|_r
\end{align*}
holds for all sequences $a,b\in\ell_r$, $G\in\ell_{\tfrac{r'}{2},\infty}$.
\end{lemnr}
\begin{bew}
	At first, Hölder's inequality gives
\begin{align*}
	\left| \sum_{k,l} a_k G_{f(k)-l} b_l \right| \leq \|a\|_r \|(G\ast b)(f(\cdot))\|_{r'}.
\end{align*}
We further observe
\begin{align*}
	\sum\limits_k |G\ast b|^{r'}(f(k)) \underset{\eqref{fastinjektiv}}{\leq}
	 C \sum\limits_l |G\ast b|^{r'}(l).
\end{align*}
Since $r<2$, we have $r'>2$, i.e. $\frac{r'}{2}>1$. Hence we apply Young's inequality with parameters $1+\frac{1}{r'}=\frac{2}{r'}+\frac{1}{r}$ to obtain
\begin{align*}
	\|G\ast b\|_{r'} \leq \|G\|_{\tfrac{r'}{2},\infty} \|b\|_r.
\end{align*}
Altogether, this provides the desired estimate. 
\end{bew}


\begin{lemnr}[Hölder's inequality in Lorentz spaces]\label{lohoe}
Let $\durch{p_i}=\durch{q_i}+\durch{r_i}$, $i=1,2$. Then
\begin{align*}
	\|fg\|_{p_1,p_2}  \lesssim \|f\|_{q_1,q_2} \|g\|_{r_1,r_2}.
\end{align*}
\end{lemnr}
\begin{bew} This is a consequence of the classical Hölder inequality and the fact that the decreasing rearrangement satisfies $(fg)^*(2t)\leq f^*(t)g^*(t)$ (see Proposition 1.4.5 No.(7) in [G]).
\end{bew}

\begin{lemnr}\label{überlapplemma}
Let $X$ be a normed vector space and $U_1,U_2,V_1,V_2\subset X$ with the property
$(U_1+V_1)\cap (U_2+V_2) \neq\emptyset$. Furthermore let $P$ be a linear functional on $X$ and let
$x_i\in U_i$, $a=P(x_2)-P(x_1)$ and $y_i\in V_i$, $b=P(y_2)-P(y_1)$.
Then we have
\begin{align}
	a+b \leq \sum_{i=1,2} [\diam(P(U_i))+\diam(P(V_i))].
\end{align}
\end{lemnr}
\begin{bew}
	Choose some $\xi\in (U_1+V_1)\cap (U_2+V_2)$. There exists	$u_i\in U_i$, $v_i\in V_i$, $i=1,2$ with
	$u_1+v_1=\xi=u_2+v_2$. It follows
\begin{align*}
	a+b=& P(x_2-x_1)+P(y_2-y_1) \\ 
	=& P(x_2-x_1)+P(y_2-y_1) +P(u_1+v_1)-P(u_2+v_2) \\
	=& P(x_2-u_2)+P(y_2-v_2)+P(u_1-x_1)+P(v_1-y_1) \\
	\leq& \diam(P(U_2))+\diam(P(V_2))+\diam(P(U_1))+\diam(P(V_1)),
\end{align*}
completing the proof.
\end{bew}

The following lemma gives us at hand a method to compare distances between two different sequences of points on the real line.

\begin{lemnr}\label{folgenlemma}
Let $I=\{0,\ldots,n\},J=\{0,\ldots,m\}\subset\N$ and let $a=(a_i)_{i\in I}\in \R^I$, $b=(b_j)_{j\in J}\in \R^J$ be two increasing, finite sequences such that
\begin{align}\label{aequidist}
	\begin{array}{ccc} C\inv \leq & a_{i+1}-a_i \leq C & \text{ for all }i\in I \\
										 C\inv \leq & b_{j+1}-b_j \leq C & \text{ for all }j\in J .\end{array}
\end{align}
Moreover we assume that 
\begin{align}\label{breite}
	b_m-b_0\leq a_n-a_0.
\end{align}

 Then there exists a function $f:J\to I$ such that
\begin{align*}
	(i)\ \quad& |a_i-b_j|\geq \durch{2} |a_i-a_{f(j)}| \qquad\\
	(ii) \quad& \forall i\in I: |\{j:f(j)=i\}|\leq 4C^2+2. \qquad
\end{align*}
\end{lemnr}
The trivial case $I=J$, $a=b$ can of course be solved by $f=id$. 
In the general setting, we in some sense had to replace every $b_j$ by some $a_{f(j)}$, leaving the distances to other points (almost) unchanged. Condition (ii) can be read as a weakening of injectivity.
\begin{bew}~\\
Without loss of generality, we may assume 
\begin{align}\label{17nov1552}
	b_m\leq a_n.
\end{align}
To be more precise, if $b_m>a_n$ would hold, we would consider $\bar a_i=-a_{n-i}$ and $\bar b_j=-b_{m-j}$, 
which also fulfill the requirements of the lemma. Then 
$$\bar b_m=-b_0 \underset{\eqref{breite}}{\leq}a_n-a_0-b_m < -a_0 =\bar a_n.$$
So, if we would find a function $\bar f$ appropriate to $\bar a$, $\bar b$ in the sense of (i) and (ii), $f(j):=n-\bar f(m-j)$ would be a solution appropriate to $a$ and $b$ since
\begin{align*}
	|a_i-b_j|=|\bar a_{n-i} -\bar b_{m-j}| \geq \durch{2}|\bar a_{n-i}-\bar a_{\bar f(m-j)}|
	=\durch{2}|a_{i}-a_{n-\bar f(m-j)}|.
\end{align*}
Furthermore, we may also assume without loss of generality
\begin{align}\label{17nov1553}
	a_0\leq b_m.
\end{align}

In the case $a_0> b_m$, we would introduce $\bar a_i=a_i-a_0+b_m < a_i$, such that $\bar a_0= b_m \leq b_m$. If $f$ is a function associated to $\bar a$ and $b$ as required, then
\begin{align*}
	|b_j-a_i|\underset{b_m<a_0}{=}a_i-b_j>\bar a_i-b_j \underset{b_m=\bar a_0}{=} |b_j-\bar a_i|
	\geq \frac{1}{2}|\bar a_{f(j)}-\bar a_i| = \frac{1}{2}|a_{f(j)}- a_i|.
\end{align*}

It would be natural for the construction of $f$ to assign an $i$ to every given $j$ in a way that $b_j$ is close to $a_i$. 
Nevertheless, if the sequences are somehow shifted against each others (for instance, $b_0\ll a_0$), there might be no $a_i$ "close" to $b_j$.
If we would always choose the closest $a_i$, this would hurt condition (ii). \\
Therefore we reflect ${a_i}'s$ at $a_0$, to ensure that for every ${b_j}$ there is a (maybe reflected) point $a_i$ nearby.\\

\stepcounter{abb}
\def\JPicScale{1}
\psset{unit=\JPicScale mm}
\psset{linewidth=0.3,dotsep=1,hatchwidth=0.3,hatchsep=1.5,shadowsize=1,dimen=middle}
\psset{dotsize=0.7 2.5,dotscale=1 1,fillcolor=black}
\psset{arrowsize=1 2,arrowlength=1,arrowinset=0.25,tbarsize=0.7 5,bracketlength=0.15,rbracketlength=0.15}
\begin{pspicture}(0,0)(150,43.75)
\psline(30,20)(100,20)
\psline(80,40)(150,40)
\psline(30,21.25)(30,18.75)
\psline(40,18.75)(40,21.25)
\psline(50,21.25)(50,18.75)
\psline(60,21.25)(60,18.75)
\psline(70,21.25)(70,18.75)
\psline(80,21.25)(80,18.75)
\psline(90,21.25)(90,18.75)
\psline(100,21.25)(100,18.75)
\psline(35,21.25)(35,18.75)
\psline(45,18.75)(45,21.25)
\psline(55,21.25)(55,18.75)
\psline(65,21.25)(65,18.75)
\psline(75,21.25)(75,18.75)
\psline(85,21.25)(85,18.75)
\psline(95,21.25)(95,18.75)
\psline(80,41.25)(80,38.75)
\psline(90,38.75)(90,41.25)
\psline(100,41.25)(100,38.75)
\psline(110,41.25)(110,38.75)
\psline(120,41.25)(120,38.75)
\psline(130,41.25)(130,38.75)
\psline(140,41.25)(140,38.75)
\psline(150,41.25)(150,38.75)
\psline[linestyle=dotted](10,40)(80,40)
\psline(10,41.25)(10,38.75)
\psline(20,38.75)(20,41.25)
\psline(30,41.25)(30,38.75)
\psline(40,41.25)(40,38.75)
\psline(50,41.25)(50,38.75)
\psline(60,41.25)(60,38.75)
\psline(70,41.25)(70,38.75)
\rput(80,36.25){$a_0$}
\rput(150,36.25){$a_n$}
\rput(103.75,18.75){$b_m$}
\rput(30,16.25){$b_0$}
\rput(88.75,16.25){$b_j$}
\rput(96.25,16.25){$b_{j+1}$}
\rput(90,43.75){$a_{f(j)}$}
\psline{->}(90,22.5)(90,36.25)
\psline{->}(95,22.5)(91.25,36.25)
\rput(45,16.25){$b_{j'}$}
\rput(40,43.75){$a_{g(j')}$}
\rput(10,36.25){$a_{-n}$}
\rput(120,43.75){$a_{f(j')}$}
\psline{->}(45,22.5)(41.25,36.25)
\rput(40,5){\textbf{Figure \thesection.\theabb:} Setting of the sequences a and b}
\end{pspicture}

Let $\bar I=I\cup(-I)$. Define the continuation of $a$ on $\bar I$ by $a_{-i}=2a_0-a_i$, $i\in I$.
We then find a function $g:J\to \bar I$ such that $|a_{g(j)}-b_j|=\min\limits_{i\in\bar I} |a_i-b_j|$. The claim is that $f=|g|$ is a solution of our problem.\\
\underline{Checking condition (i):}\\
{Case 1: $g(j)\geq0$.} Here
\begin{align*}
	|a_i-a_{f(j)}|=|a_i-a_{g(j)}|\leq|a_i-b_j|+|b_j-a_{g(j)}| 
				\underset{\text{minimality of }g}{\leq}2|a_i-b_j|.
\end{align*}
{Case 2: $g(j)<0$.} Here
\begin{align*}
	|a_i-a_{f(j)}|=|a_i-a_{-g(j)}|=&|a_i-2a_0+a_{g(j)}| \\
	\leq& |a_i-a_0|+|a_0-a_{g(j)}|
	\underset{\text{monotony}}{=} a_i-a_0+a_0-a_{g(j)} = |a_i-a_{g(j)}|,
\end{align*}
and we proceed as in case 1.\\
\underline{Checking condition (ii):}\\
Obviously, it is sufficient to show that $|\{j:g(j)=i\}|\leq 2C^2+1$ holds for every $i\in \bar I$.
Thus let $g(j)=i$.\\
We claim that $|a_i-b_j|\leq C$ and check this:\\
{Case 1: $a_i<b_j$.} If $i=n$, then 
$b_j \leq b_m \underset{\eqref{17nov1552}}{\leq} a_n = a_i<b_j$, hence $i<n$. Thus $a_{i+1}$ is well-defined and $b_j<a_{i+1}$ holds according to the minimality in the choice of $g$. 
We conclude $|a_i-b_j|=b_j-a_i\leq a_{i+1}-a_i \leq C$.\\
{Case 2: $a_i>b_j$.} Would $i=-n$, then
$$b_j<a_i=a_{-n}=2a_0-a_n\underset{\eqref{17nov1553}}{\leq}b_m-(a_n-a_0)\underset{\eqref{breite}}{\leq}b_0\leq b_j.$$ Hence we have $i>-n$, thus $a_{i-1}$ is well-defined and $b_j>a_{i-1}$ holds. 
We finish as in case 1.\\
\underline{Case 3: $a_i=b_j$.} This case is trivial.\\
If additionally $g(j')=i$, then $|a_i-b_{j}|\leq C$ and $|a_i-b_{j'}|\leq C$.
Now we utilise \eqref{aequidist}. Therefrom we get by induction
\begin{align*}
	\durch{C}|j-j'|\leq|b_j-b_{j'}|\leq|b_j-a_i|+|a_i-b_{j'}|\leq 2C,
\end{align*}
i.e.
\begin{align*}
	|j-j'|\leq 2C^2,
\end{align*}
and thus condition (ii).
\end{bew}

\begin{bem}
Of course, if assumption \eqref{breite} is not fulfilled, we may apply the lemma with interchanged roles of $I$ and $J$ or $a$ and $b$ respectively. Anyway, the lemma actually remains essentially valid if these quantities are \emph{not} interchanged. For this, we just need some further reflections of the $a_i$-sequence, if necessary also at the "upper" end at $a_n$. This would enlarge the bound $4C^2+2$ from (ii) depending on the relation between $n$ and $m$. 
However, with to many reflections, the notation would quickly become confusing, thus I omit such a version of the lemma.
\end{bem}

\end{document}